\documentclass[a4paper, reqno, 10pt]{report}

\usepackage[sc, noBBpl]{mathpazo}	
\usepackage[bbgreekl]{mathbbol}		
\usepackage{amssymb, amsmath, amsthm}

\usepackage[T1]{fontenc}
\usepackage[utf8]{inputenc}
\usepackage{graphicx}
\usepackage{comment}
\usepackage[linkbordercolor={0.8 0.9 0.9}, citebordercolor={0.85 0.85 0.85}]{hyperref}
\usepackage[normalem]{ulem}

\usepackage{amsmath, amsthm, amssymb} 
\usepackage[all]{xy}
\usepackage[usenames,dvipsnames]{color}
\usepackage{bibentry}
\nobibliography*

\numberwithin{equation}{section}

\theoremstyle{plain}
\newtheorem{theorem}{Theorem}[section]
\newtheorem{proposition}[theorem]{Proposition}
\newtheorem{lemma}[theorem]{Lemma}
\newtheorem{corollary}[theorem]{Corollary}
\theoremstyle{definition}
\newtheorem{definition}[theorem]{Definition}
\newtheorem{notation}[theorem]{Notation}
\newtheorem{example}[theorem]{Example}
\theoremstyle{remark}
\newtheorem{remark}[theorem]{Remark}

\newcommand{\st}{\;|\;}

\newcommand{\half}{{\textstyle\frac{1}{2}}}
\newcommand{\ip}[1]{\langle #1 \rangle}

\newcommand{\into}{\hookrightarrow}
\newcommand{\onto}{\twoheadrightarrow}

\newcommand{\Thus}{\quad\Longrightarrow\quad}

\newcommand{\lie}{\mathfrak}

\renewcommand{\AA}{\mathbb{A}}
\newcommand{\CC}{\mathbb{C}}

\newcommand{\JJ}{\mathbb{J}}

\newcommand{\KK}{\mathbb{K}}
\newcommand{\LL}{\mathbb{L}}
\newcommand{\NN}{\mathbb{N}}

\newcommand{\RR}{\mathbb{R}}

\newcommand{\XX}{\mathbb{X}}
\newcommand{\ZZ}{\mathbb{Z}}

\newcommand{\sA}{\mathcal{A}}
\newcommand{\sB}{\mathcal{B}}
\newcommand{\sC}{\mathcal{C}}
\newcommand{\sD}{\mathcal{D}}
\newcommand{\sE}{\mathcal{E}}
\newcommand{\sF}{\mathcal{F}}

\newcommand{\sS}{\mathcal{S}}
\newcommand{\sU}{\mathcal{U}}
\newcommand{\sX}{\mathcal{X}}

\newcommand{\bfa}{\mathbf{a}}
\newcommand{\bfb}{\mathbf{b}}

\newcommand{\rmA}{\mathrm{A}}
\newcommand{\rmB}{\mathrm{B}}
\newcommand{\rmC}{\mathrm{C}}
\newcommand{\rmG}{\mathrm{G}}

\newcommand{\SL}{\mathrm{SL}}
\newcommand{\SU}{\mathrm{SU}}

\newcommand{\Sp}{\mathrm{Sp}}

\newcommand{\bP}{\mathbf{P}}		
\newcommand{\bQ}{\mathbf{Q}}		
\newcommand{\bDelta}{\mathbf{\Delta}}	
\newcommand{\bSigma}{\mathbf{\Sigma}}	

\newcommand{\Irr}{\mathrm{Irr}}

\renewcommand{\Re}{\mathop{\mathrm{Re}}\nolimits}

\DeclareMathOperator{\ph}{ph}
\DeclareMathOperator{\supp}{supp}
\DeclareMathOperator{\End}{End}

\DeclareMathOperator{\Ind}{Ind}
\DeclareMathOperator{\res}{res}

\DeclareMathOperator{\id}{id}

\DeclareMathOperator{\ad}{ad}

\DeclareMathOperator{\Op}{Op}
\DeclareMathOperator{\gr}{gr}
\DeclareMathOperator{\singsupp}{sing-supp}

\newcommand{\pp}{\mathbb{p}}
\newcommand{\qq}{\mathbb{q}}

\newcommand{\pr}{\mathrm{pr}}
\newcommand{\DO}{\mathrm{DO}}
\newcommand{\CP}{\CC \mathrm{P}}

\newcommand{\Uqg}{\sU_q(\lie{g})}
\newcommand{\UqRk}{\sU_q(\lie{k})}
\newcommand{\UqRl}{\sU_q(\lie{l})}
\newcommand{\slot}{\,\cdot\,}

\newcommand{\Cc}{C_\mathrm{c}}
\newcommand{\Cp}{C_\mathrm{p}}
\newcommand{\sAc}{\sA_\mathrm{c}}
\newcommand{\sDp}{\sD_\mathrm{p}}

\newcommand{\Ell}{\mathcal{E}\!\ell\!\ell}
\newcommand{\Dom}{\mathrm{Dom}}

\newcommand{\tHM}{\lie{t}_HM}
\newcommand{\tsHM}{\lie{t}^*_HM}
\newcommand{\THM}{T_HM}
\newcommand{\ttHM}{\mathbb{t}_HM}
\newcommand{\TTHM}{\mathbb{T}_HM}
\newcommand{\PPsi}{\mathbb{\Psi}}
\newcommand{\BGG}{[\mathrm{BGG}]}
\newcommand{\bfmu}{{\boldsymbol{\mu}}}
\newcommand{\op}{\mathrm{op}}
\newcommand{\cop}{\mathrm{cop}}
\newcommand{\hit}{\!\rightharpoonup\!}

\begin{document}

\title{On pseudodifferential operators on filtered and multifiltered manifolds}
\author{Robert Yuncken}


\thispagestyle{empty}

\begin{center}
	{\Large\bf UNIVERSIT\'E CLERMONT AUVERGNE} 
	
	\medskip
	{\bf \large Laboratoire de Math\'ematiques Blaise Pascal \\ UMR 6620 CNRS}
	
	\medskip 
	
	\hrule
	\smallskip
	\hrule
	
	\bigskip\bigskip  \bigskip
	
	{\Large TRAVAUX DE RECHERCHE }
	
	\bigskip \bigskip 
	
	{\large pr\'esent\'es en vue de l'obtention de }
	
	\bigskip  \bigskip
	
	{\Large L'HABILITATION \`A DIRIGER DES RECHERCHES EN MATH\'EMATIQUES}
	
	\bigskip \bigskip
	
	{\Large\bf ON PSEUDODIFFERENTIAL OPERATORS ON FILTERED AND MULTIFILTERED MANIFOLDS}
	
	\bigskip \bigskip
	
	{\large  par}
	
	\bigskip \bigskip
	
	{\large\bf  Robert Yuncken }
	
	\smallskip
	{\sl Ma\^itre de conf\'erences en Math\'ematiques \`a l'Universit\'e Clermont Auvergne}
\end{center}

\bigskip \bigskip \noindent  Rapporteurs :

\bigskip 
\begin{itemize}
	\item Georges Skandalis, Universit\'e Paris 7 Denis Diderot
	\item Pierre Julg, Universit\'e d'Orleans
	\item Sergey Neshveyev, Univerist\'e d'Oslo
\end{itemize}

\bigskip \bigskip \noindent Habilitation soutenue publiquement le 28 septembre 2018 devant le jury compos\'e de :

\bigskip 
\begin{itemize}
	\item Georges Skandalis, Universit\'e Paris 7 Denis Diderot
  \item Ryszard Nest, University de Copenhague
  \item Pierre Julg, Universit\'e d'Orl\'eans 
  \item Claire Debord, Universit\'e Paris 7 Denis Diderot
  \item Julien Bichon, Universit\'e Clermont Auvergne
\end{itemize}

\newpage 
\thispagestyle{empty}
\vspace*{\fill}

\noindent {\large  
	
	\bigskip \noindent
	UNIVERSIT\'E CLERMONT AUVERGNE \\
	Laboratoire de Math\'ematiques Blaise Pascal \\
	UMR 6620 CNRS \\
	Campus des Cézeaux \\
	63178 Aubière cedex \\ 
	France}

\bigskip \noindent 
T\'el\'ephone : +33 (0)4 73 40 76 97 \\
Courriel : \url{robert.yuncken@uca.fr}

\begin{abstract}
 We present various different approaches to constructing algebras of pseudodifferential operators adapted to particular geometric situations.  A general goal is the study of index problems in situations where standard elliptic theory is insufficient.  We also present some applications of these constructions.

 We begin by presenting a characterization of pseudodifferential operators in terms of distributions on the tangent groupoid which are essentially homogeneous with respect to the natural $\RR^\times_+$-action.  This is carried out in the generality of filtered manifolds, which are manifolds that are locally modelled on nilpotent groups, generalizing contact, CR and parabolic geometries.  We describe the tangent groupoid of a filtered manifold, and use this to construct a pseudodifferential calculus analogous to the unpublished calculus of Melin.
 
 Next, we describe a rudimentary multifiltered pseudodifferential theory on the full flag manifold $\sX$ of a complex semisimple Lie group $G$ which allows us to simultaneously treat longitudinal pseudodifferential operators along every one of the canonical fibrations of $\sX$ over smaller flag manifolds.  The motivating application is the construction of a $G$-equivariant $K$-homology class from the Bernstein-Gelfand-Gelfand complex of a semisimple group.  This construction been completely resolved for only a few groups, and we will discuss the remaining obstacles as well as the successes.
 
 Finally, we discuss pseudodifferential operators on two classes of quantum flag manifolds.   First, we consider quantum projective spaces $\CP^n_q$, where we can generalize the abstract pseudodifferential theory of Connes and Moscovici to obtain a twisted algebra of pseudodifferential operators associated to the Dolbeault-Dirac operator.  Secondly, we look at the full flag manifolds of $\SU_q(n)$, where we instead generalize the multifiltered construction of the classical flag manifolds, thus obtaining an equivariant fundamental class for the full flag variety of $\SU_q(3)$ from the Bernstein-Gelfand-Gelfand complex.  As applications, we obtain equivariant Poincar\'e duality of the quantum flag manifold of $\SU_q(3)$ and the Baum-Connes Conjecture for the discrete dual of $\SU_q(3)$ with trivial coefficients.

\end{abstract}

\shipout\null

{\fontfamily{lmss}\selectfont 
\renewcommand{\abstractname}{Thanks}
\begin{abstract}
 
There are a lot of people to thank.  
Let me start at the beginning by again thanking my PhD adviser, Nigel Higson, for his guidance and good taste.  He has been a constant source of inspiration.  Thanks also to the late great John Roe, another big inspiration, who  has had far more influence on my work than may be apparent.
\medskip

Huge thanks to my collaborators---Heath, Christian, Erik, Julien, Marco and Robin---who have all made my mathematics better, and more fun.
\medskip

Thanks to the French operator algebra community for adopting me.  I couldn't ask for a nicer clan.  The list is long, but let me particularly thank Herv\'e, Jean-Marie, Claire, Georges, Saad, Jean and Julien.
\medskip

Thanks very much to the referees and jury members who are offering their precious time for assessing this work.
\medskip

And last but not least, thanks to my family: to my parents for their support, to my kids for keeping me happy with their bottomless spirit, and  to HJ,~who owns a piece of everything I have done and who deserves more thanks than I could fit on this page.

\end{abstract}
}

\shipout\null

\thispagestyle{empty}
\begin{center}
{\Large\bf List of publications presented}
\end{center}

\vspace{8mm}

The following are articles which will be presented in this memoir.  

\subsection*{Chapter \ref{sec:groupoid_PsiDOs}: Pseudodifferential operators from tangent groupoids}
\begin{description}
	\item[\cite{VanYun:groupoid}] \bibentry{VanYun:groupoid}
	\item[\cite{VanYun:PsiDOs}] \bibentry{VanYun:PsiDOs}
\end{description}

\subsection*{Chapter \ref{sec:flags}: Pseudodifferential operators on multifiltered manifolds and equivariant index theory for complex semisimple groups}
\begin{description}
	\item[\cite{Yuncken:foliations}] \bibentry{Yuncken:foliations}
\end{description}
This work is heavily motivated by the earlier publication:
\begin{description}
	\item[\cite{Yuncken:BGG}] \bibentry{Yuncken:BGG}
\end{description}

\subsection*{Chapter \ref{sec:CQG}: Pseudodifferential operators on quantum flag manifolds}
\begin{description}
	\item[\cite{MatYun:PsiDOs}] \bibentry{MatYun:PsiDOs}
	\item[\cite{VoiYun:SUq3}] \bibentry{VoiYun:SUq3}
\end{description}
Material from the following preprint will also be used:
\begin{description}
	\item[\cite{VoiYun:CQG}] \bibentry{VoiYun:CQG}
\end{description}

\newpage

\shipout\null

\setcounter{page}{1}
\tableofcontents


\chapter{Introduction}

The goal of this memoir is to describe a handful of different approaches to constructing algebras of pseudodifferential operators on manifolds and noncommutative spaces with particular geometries.  In each case, the methods will be rather different, but the goals similar.   We will also present some of the applications of these pseudodifferential algebras.

Pseudodifferential operators are, of course, a tool for studying differential operators, particularly elliptic differential operators and their generalizations.   Our motivations come from index theoretic problems where the underlying geometry demands more exotic pseudodifferential theories than ordinary elliptic theory.
We will be particularly interested in:
\begin{enumerate}
	\item[(1)] Rockland operators, including subelliptic operators on contact and\linebreak Heisenberg manifolds;
	\item[(2)] the Bernstein-Gelfand-Gelfand complex, a canonical complex of equivariant differential operators on the flag manifold of a semisimple Lie group;
	\item[(3)] analogues of the above operators on quantum homogeneous spaces.
\end{enumerate}

As mentioned, the techniques of analysis for the various examples will be very different: in (1) we shall use groupoid techniques, in (2) we use convolution operators on nilpotent Lie groups, and in (3) methods of noncommutative geometry and noncommutative harmonic analysis.  But in each case, the general strategy---the requirements of the associated pseudodifferential calculus---are roughly similar.

We shall dedicate this introduction to a very broad discussion of this strategy.    Our goal here is to isolate the general properties of pseudodifferential operators which are essential to index theory.  Finally, at the end of the introduction, we will give a brief overview of the particular problems to be discussed in the following chapters.

\medskip


\section{Elliptic operators and their generalizations}

Index theory begins with the study of analytic properties of elliptic differential operators on a closed manifold $M$.  The crucial properties of elliptic operators, which for the moment we will state imprecisely%
\footnote{For a precise interpretation of these statements, one should consider an elliptic operator of order $m$ as a bounded operator between Sobolev spaces $H^{s+m}(M) \to H^s(M)$ for any $s\in\RR$.}
, are:
\begin{enumerate}
 \item[(1)] Elliptic differential operators on a closed manifold are Fredholm.

 \item[(2)] If two elliptic differential operators have the same principal symbol then their difference is a compact operator.
\end{enumerate}
These two properties imply that the Fredholm index of an elliptic differential operator depends only upon the principal symbol of the operator.  Moreover, the index is unchanged by smooth perturbations of the principal symbol.  In other words, the index of an elliptic operator depends only on some topological data associated to the principal symbol.
The appropriate topological data is the $K$-theory class of the principal symbol $\sigma(D)$, and for elliptic operators, the index is calculated by the famous Atiyah-Singer Formula:
\[
 \Ind(D) = \int_M \mathrm{ch}(\sigma(D)) \mathrm{Td}(M) .
\]

Properties (1) and (2) are valid for a much larger class of operators than elliptic differential operators.  The examples we consider will all be \emph{Rockland} (see Definition \ref{def:Rockland}).  This is a generalization of ellipticity to filtered manifolds, \emph{i.e.}\ manifolds with a filtration of the tangent bundle compatible with the Lie bracket of vector fields.  Such operators appear naturally in applications, particularly in CR and contact geometry, and more recently in parabolic geometry, see \emph{e.g.} \cite{FolSte:estimates,Rockland,BeaGre,Taylor:microlocal, Melin:preprint, EpsMelMen, CapSloSou:BGG, DavHal:BGG}.

One example of particular interest to us is the Bernstein-Gelfand-Gelfand complex.  This is a canonical equivariant differential complex on a flag manifold of a semisimple Lie group, or more generally on a parabolic manifold \cite{CapSloSou:BGG}.  The Bernstein-Gelfand-Gelfand complex has gained much attention recently, originally due to its appearance in twistor theory \cite{BasEas}.  From our point of view, the BGG complex is an obligatory replacement of a Dirac-type operator for the equivariant index theory of semisimple Lie groups and their quantizations; see Sections \ref{sec:flags} and \ref{sec:CQG}.

\medskip

Another important property of elliptic differential operators, also shared by Rockland operators, is hypoellipticity.
Let $E$ be a vector bundle over a manifold $M$ without boundary.  We write $C^\infty(M;E)$ for the space of smooth sections of $E$, and $\sD'(M;E)$ for the space of distributional sections.

\begin{definition}
  A linear differential operator $P:\sD'(M;E)\to \sD'(M;E)$ is \emph{hypoelliptic} if, for any open set $U\subseteq M$ and any distribution $u\in\sD'(M;E)$ we have
	\[
	  Pu|_U \in C^\infty(U;E) \Thus u|_U \in C^\infty(U;E).
	\]
\end{definition}

In other words, hypoellipticity guarantees smooth solutions $u$ to the partial differential equation $Pu=f$ whenever the right-hand-side $f$ is smooth.  Both hypoellipticity and Fredholmness are typically proven by constructing a pseudodifferential calculus adapted to $P$.

\section{Pseudodifferential operators}
\label{sec:calculi}

Let $M$ be a smooth manifold without boundary.
For simplicity, we will suppress coefficient vector bundles in this section.  We will also be lazy about the difference between functions and densities.   We will be more precise in later chapters.

\begin{definition}
 A distribution $p \in \sD'(M \times M)$ is called
 \begin{itemize}
  \item \emph{properly supported} if the restriction to $\supp(p)$ of each of the two coordinate projections $\pr_1, \pr_2: M \times M \onto M$ is a proper map;
  \item \emph{semiregular in both variables} if $u\in C^\infty(M) \bar\otimes \sD'(M) \;\cap\; \sD'(M) \bar\otimes C^\infty(M)$.
 \end{itemize}
 We write $\sDp'(M\times M)$ for the space of distributions on $M\times M$ which are properly supported and semiregular in both variables---for a better context see  Definition \ref{def:proper}.
\end{definition}

Via the Schwartz kernel theorem, $\sDp'(M\times M)$ is in bijection with the algebra of operators on $\sD'(M)$ which preserve each of the subspaces $\sE'(M)$, $C^\infty(M)$ and $\Cc^\infty(M)$, see \cite{Treves:TVS}.  By abuse of notation, we will often identify such an operator with its kernel in $\sDp'(M\times M)$.  The product of such operators corresponds to the convolution product of distributions:
\[
 p*q(x,z) = \int_M p(x,y)q(y,z) \,dy, \qquad p,q \in \sD'_p(M \times M),
\]
which makes sense thanks to semiregularity.
In particular, linear differential operators on $M$ correspond to elements of $\sDp'(M \times M)$ with support on the diagonal.

Now let $\Ell$ be some set of linear differential operators on $M$.  We are imagining a class of operators for which we hope to prove hypoellipticity or Fredholmness.  The following definition is intended only as a guiding philosophy.


\begin{definition}
 \label{def:PsiDOs}
 An \emph{algebra of pseudodifferential operators adapted to $\Ell$} will mean a $\ZZ$-filtered algebra $\Psi^\bullet \subset \sDp'(M\times M)$ which has the following properties:
 \begin{enumerate}
  \item $\Psi^\bullet$ contains (Schwartz kernels of) all linear differential operators on $M$.
  \item Pseudolocality: Every $p\in\Psi^\bullet$ is equal to a smooth function off the diagonal.
  \item $\Psi^{-\infty} = \bigcap_{m\in\ZZ} \Psi^m  = \Cp^\infty(M \times M)$ is the algebra of properly supported smoothing kernels.
  \item Existence of asymptotic limits: Let $m\in\ZZ$ and $p_i \in \Psi^{m-i}$ for all $i\in\NN$  Then there exists $p \in \Psi^m$ such that for every $n\in\ZZ$,
  \[
   p - \sum_{i=0}^k p_i \in \Psi^{-n} \qquad \text{ for all } k \gg 0.
  \]
  \item Existence of parametrices: If $p \in \Ell$, then there exists $q\in \Psi^{\bullet}$ such that 
  \[
   I - p*q, ~ I - q*p ~ \in \Psi^{-1},
  \]
  where $I$ is the Schwartz kernel of the identity operator.
 \end{enumerate}
\end{definition}

For instance, the algebra of classical (step one polyhomogeneous) pseudodifferential operators on $M$ is adapted to the class of elliptic differential operators in this sense.  If $M$ is a contact manifold then the Heisenberg calculus of Beals and Greiner \cite{BeaGre} is adapted to the maximally hypoelliptic (\emph{i.e.}, Rockland) operators.

The point of the five properties above is that they imply the hypoellipticity of operators in $\Ell$.

\begin{theorem}
\label{thm:Ell}
 Let $\Ell$ be a class of linear differential operators on a manifold $M$ and suppose that there exists an algebra of pseudodifferential operators $\Psi^\bullet$ adapted to $\Ell$ in the sense of Definition \ref{def:PsiDOs}.  Then the operators in $\Ell$ are hypoelliptic.
\end{theorem}

\begin{proof}
 Let $P\in\Ell$ of order $m$.  Let $Q \in \Psi^{-m}$ be a parametrix, and put $R=I-PQ \in \Psi^{-1}$.  The series $Q\sum_{k=0}^\infty R^k$ admits an asymptotic limit, which we denote by $\tilde{Q}$.  One checks that $I-P \tilde{Q}$, $I-\tilde{Q}P \in \Psi^{-\infty}$. 

 Now, for any $u\in \sD'(M)$ we have
 \[
  u = (I-\tilde{Q}P)u + Pu,
 \]
 and it follows that $\singsupp(u) \subseteq \singsupp(Pu)$.  Thus $P$ is hypoelliptic.
\end{proof}

A similar philosophy could be applied to Fredholmness.
By adding slightly more structure to Definition \ref{def:PsiDOs}, replicating basic elements of Sobolev theory, we could deduce Fredholmness of the operators in $\Ell$.  We will not formulate a precise statement here.

\section{Overview}

The goal of the following chapters, roughly speaking, will be to obtain some version of the properties stated in Definition \ref{def:PsiDOs} for algebras of pseudodifferential operators which are adapted to particular differential operators.  This will be achieved with more or less success depending on the examples. 

\medskip

In Chapter \ref{sec:groupoid_PsiDOs}, we will describe an approach to pseudodifferential theory based on the tangent groupoid of a filtered manifold \cite{VanYun:PsiDOs}.  The main result is a simple groupoid definition of pseudodifferential kernels which reproduces the classical calculus, the Heisenberg calculus \cite{BeaGre, Taylor:microlocal}, and more generally, Melin's unpublished pseudodifferential calculus on an arbitrary filtered manifold \cite{Melin:preprint}.  Philosophically, this result is supposed to indicate a general principal: in order to define a pseudodifferential calculus, it is sufficient to construct a tangent groupoid adapted to the geometry, from which the pseudodifferential calculus follows automatically.

\medskip

In Chapter \ref{sec:flags}, we consider longitudinal pseudodifferential operators on manifolds equipped with multiple foliations \cite{Yuncken:foliations}.  The examples of interest are the flag manifolds of complex semisimple Lie groups, which admit a family of canonical fibrations over smaller flag manifolds. The results presented are designed to simultaneously treat longitudinally pseudodifferential operators along the fibres of different fibrations.

The theory we present is far from being a pseudodifferential calculus in the sense of Definition \ref{def:PsiDOs}, since it essentially distinguishes only between longitudinal pseudodifferential operators of order $0$ and those of negative order.  Still, it  allows us to interpret the Bernstein-Gelfand-Gelfand complex as an equivariant fundamental class in the Kasparov $K$-homology of the full flag manifold of $\SL(3,\CC)$.  As has been shown in \cite{Yuncken:BGG}, this leads an explicit construction of the $\gamma$ element from the Baum-Connes Conjecture.  

\medskip

Finally, In Chapter \ref{sec:CQG}, we pass to pseudodifferential theory in noncommutative geometry---specifically, on quantized flag manifolds.  These noncommutative spaces would seem to deserve the status of ``noncommutative manifolds'', but it has been very difficult to incorporate them into Connes' framework.  Much of the difficulty is due to the absence of a reasonable notion of pseudodifferential calculus for these spaces.  

We will consider two classes of quantum flag manifolds, with rather different behaviours.  Firstly, we look at the quantized projective spaces $\CP^n_q$, where things are relatively simple.  We show that the existing spectral triples on $\CP^n_q$ \cite{Krahmer:Dirac,KraTuc,DAnDab} are twisted regular, in the sense of \cite{ConMos:twisted}.  This is achieved by means of a generalization of the abstract pseudodifferential calculus of Connes and Moscovici \cite{ConMos:local_index_formula, ConMos:twisted, Moscovici:twisted}.

Secondly, we consider the full flag manifold $\sX_q$ of $\SU_q(n)$, for which we have no candidates for a spectral triple.  Instead, we show that essential elements of the multi-filtered harmonic analysis of Chapter \ref{sec:flags} can be carried over to quantum flag manifolds \cite{VoiYun:SUq3}.  As a result, the quantum analogue of the Bernstein-Gelfand-Gelfand complex for $\SU_q(3)$ yields a fundamental class in equivariant $K$-homology $K^{\SL_q(3,\CC)}(\sX_q,\CC)$, which allows us to prove equivariant Poincar\'e Duality for $\sX_q$ and the Baum-Connes Conjecture (with trivial coefficients) for the discrete dual of $\SU_q(3)$.

\bigskip

This document is intended as a survey.  Proofs, when given, will be sketched, with references to the original publications for details.

\chapter{Pseudodifferential operators from tangent groupoids}
\label{sec:groupoid_PsiDOs}

Lie groupoids are a fusion of differential geometry and algebra.  This structure allow us to describe linear operators on smooth manifolds as convolution operators.  In particular, pseudodifferential operators on a smooth manifold $M$ can be interpreted as convolution operators by their Schwartz kernels on the pair groupoid $M\times M$, while the principal symbol, up to Fourier transform, is a convolution operator on the tangent bundle $TM$  seen as a bundle of abelian Lie groups.  Connes introduced the tangent groupoid \cite{Connes:NCG} as a tool for smoothly deforming a pseudodifferential operator to its principal symbol. 

Debord and Skandalis realized that classical pseudodifferential operators could be characterized in terms of the tangent groupoid \cite{DebSka:Rx-action}.
Inspired by this, we showed in  \cite{VanYun:PsiDOs} that tangent bundles can be used to define pseudodifferential calculi adapted to a wide variety of geometric situations.  
The context we work in is that of a \emph{filtered manifold} $M$, which is a structure generalizing contact manifolds, Heisenberg manifolds, and parabolic manifolds.  Pseudodifferential calculi on such manifolds have been developed by many authors  \cite{FolSte:estimates, Taylor:microlocal, BeaGre, EpsMel:book, Melin:preprint}.

One can construct a tangent groupoid $\TTHM$ for such manifolds \cite{VanErp:thesis, Ponge:groupoid, VanYun:groupoid, ChoPon}.  It is a deformation of the pair groupoid $M\times M$  to the bundle of osculating groups of $M$.  As in \cite{DebSka:Rx-action}, it admits a one-parameter group of automorphisms $(\alpha_\lambda)_{\lambda\in \RR^\times_+}$.  We define a pseudodifferential kernel of order $m$ to be a properly supported semiregular Schwartz kernel $p$ which is the restriction to $M\times M$ of a smooth family $\pp$ of distributions on the $r$-fibres of $\TTHM$ satisfying
\[
 \alpha_{\lambda*}\pp - \lambda^m\pp \in C^\infty(\TTHM), \qquad \forall \lambda \in \RR^\times_+.
\]
A precise statement will be given in Section \ref{sec:H-PsiDOs}.
We showed in \cite{VanYun:PsiDOs} that this coincides with the usual definitions for the classical and Heisenberg calculi.  This chapter is dedicated to explaining this definition and its consequences.

\section{Filtered manifolds}
\label{sec:filtered_manifolds}

\begin{definition}
 A \emph{filtered manifold} (also called a \emph{Carnot manifold} in \cite{ChoPon}) is a $C^\infty$ manifold $M$ equipped with a filtration of its tangent bundle by subbundles
 \[
  \mathbf{0} = H^0 \leq H^1 \leq \cdots \leq H^N = TM,
 \]
 such that the space of smooth vector fields $\Gamma^\infty(TM)$ becomes a filtered Lie algebra:
 \[
  [\Gamma^\infty(H^i), \Gamma^\infty(H^j)] \subseteq \Gamma^\infty(H^{i+j}) \qquad\qquad \forall i,j.
 \]
 Here we are using the convention $H^k = TM$ for all $k\geq N$.

 The number $N$ will be called the \emph{depth} of the filtration.  Smooth sections of $H^k$ are called \emph{vector fields of $H$-order $\leq k$}.
\end{definition}

The filtration of vector fields by $H$-order generates an algebra filtration on the differential operators\footnote{If $M$ is noncompact, we will insist that $\DO(M)$ consists of the differential operators of finite order.}  $\DO(M)$.   The set of differential operators of $H$-order $\leq m$ will be denoted $\DO_H^m(M)$.

\begin{example}
 Any manifold $M$ can be equipped with the \emph{trivial filtration} of depth $1$, where $H^1 = TM$.  The resulting filtration on $\DO(M)$ corresponds to the usual notion of order of a differential operator.
\end{example}

\begin{example}
 Let $M$ be the Heisenberg group of dimension $3$.  We write $X,Y,Z$ for the usual left-invariant vector fields which satisfy $[X,Z]=[Y,Z]=0$ and $[X,Y]=Z$.  We obtain a filtration of $TM$ of depth $2$ by defining $H^1$ to be the $2$-dimensional subbundle spanned by $X$ and $Y$.  Then the vector fields $X$ and $Y$ are order $1$, while $Z=XY-YX$ counts as order $2$.

 This example generalizes to the notion of \emph{Heisenberg order} for differential operators on a contact manifold, see \cite{BeaGre}.
\end{example}

\begin{example}
 Any subbundle $H$ of the tangent bundle of $M$ gives rise to a filtered manifold of depth $2$,
 \[
   \mathbf{0} \leq H \leq TM,
 \]
 since the condition on Lie brackets is trivial.  This includes contact manifolds, where $H$ is the contact hyperplane bundle, as well as foliated manifolds, where $H$ is the tangent space to the leaves. 
\end{example}

\begin{example}
 Flag manifolds and, more generally, parabolic geometries are examples of filtered manifolds \cite{CapSlo}.
\end{example}

The principal part of a differential operator on a filtered manifold is defined via the passage from the filtered algebra $\DO_H(M)$ to its associated graded algebra.  We fix some general notation.

\begin{notation}
 For any $\NN$-filtered vector space (or vector bundle) $V^\bullet$, we write 
\[
 \gr V = \bigoplus_m \gr_m V := \bigoplus_m V^m/V^{m-1}
\]
for the associated graded space (or graded bundle) and 
\[
 \sigma_m : V^m \to \gr_m V \into \gr V
\]
for the canonical quotient maps.
\end{notation}

\begin{definition}
 \label{def:principal_part}
If $P\in \DO_H^m(M)$ is a differential operator of $H$-order $m$, its \emph{principal part} is defined as $\sigma_m(P)$.  
\end{definition}

As in elliptic theory, it is the principal part of a differential operator which will be used to prove hypoellipticity and Fredholmness.  The usual criterion of ellipticity needs to be replaced by the Rockland condition, which we will state in Section \ref{sec:Rockland}. 

\section{Schwartz kernels and fibred distributions on Lie groupoids}
\label{sec:tangent_groupoid}

\begin{definition}
 The \emph{pair groupoid} of a smooth manifold $M$ is the Cartesian product $M\times M$ equipped with range and source maps
 \[
  r(x,y) = x \qquad s(x,y) = y
 \]
 and groupoid laws
 \[
  (x,y)(y,z) = (x,z), \qquad (x,y)^{-1} = (y,x).
 \]
\end{definition}

This is the appropriate structure for describing Schwartz kernels of continuous operators on $C^\infty(M)$; see below.
Recall, from Section \ref{sec:calculi} that we are interested in Schwartz kernels which are properly supported and semiregular in both variables.  These notions have groupoid interpretations, which will important in the sequel.  

\begin{definition}
 Let $G$ be a Lie groupoid with base space $G^{(0)}$.
 Consider $C^\infty(G)$ as a $C^\infty(G^{(0)})$-module via the range fibration:
 \begin{equation}
  \label{eq:r-module}
  a.f  = r^*(a)f, \qquad\qquad a\in C^\infty(G^{(0)}), ~ f \in C^\infty(G).
 \end{equation}
 Then an \emph{$r$-fibred distribution}%
 \footnote{
  Strictly speaking, we are defining here the \emph{$r$-fibred distributions with compact vertical support} (or \emph{$r$-proper support}).  We will not need any more general support conditions.
 }
 on $G$ means a continuous $C^\infty(G^{(0)})$-linear operator
 \begin{equation}
  u : C^\infty(G) \to C^\infty(G^{(0)}).
 \end{equation}
 The space of $r$-fibred distributions on $G$ is denoted by $\sE'_r(G)$.

 Likewise, we define the space $\sE_s(G)$ of \emph{$s$-fibred distributions} on $G$, where we replace the $C^\infty(G^{(0)})$-module structure \eqref{eq:r-module} with its analogue using the source fibration.
\end{definition}


It is possible to integrate an $r$- or $s$-fibred distribution to an ordinary distribution on $G$ by composing with a fixed choice of nowhere zero smooth density $\mu$ on the base space $G^{(0)}$.  This yields maps
\begin{align*}
 \mu_r : &\sE_r'(G) \to \sD(G) ;  \qquad u \mapsto \mu\circ u, \\
 \mu_s : &\sE_s'(G) \to \sD(G) ;  \qquad v \mapsto \mu\circ v. 
\end{align*}
These maps depend on the choice of $\mu$, but their images do not.  In this way, we may see the $r$-fibred or $s$-fibred distributions as those distributions on $G$ which can be disintegrated as a smooth family of distributions on the $r$-fibres or $s$-fibres, respectively.  We will want to have both.

\begin{definition}
 \label{def:proper}
 A distribution $w\in\sD'(G)$ will be called \emph{proper} if it lies in the image of both $\mu_r$ and $\mu_s$.  We write $\sDp'(G)$ for the space of proper distributions.  We also define the spaces of \emph{proper $r$-fibred} and \emph{proper $s$-fibred distributions}, respectively, as 
 \begin{align*}
  \sE'_{r,s}(G) &= \mu_r^{-1}\sDp'(G) \qquad\qquad \text{and}&
  \sE'_{s,r}(G) &= \mu_s^{-1}\sDp'(G).
 \end{align*}
\end{definition}

To see the importance of these definitions, consider again the pair groupoid $G= M\times M$.  We will fix, once and for all, a nowhere vanishing smooth density $\mu$ on $M$, which allows us to identify $C^\infty(M)$ as a subspace of $\sD'(M)$.  Recall that if $p \in \sD'(M \times M)$ is an arbitrary distribution, the Schwartz kernel operator with kernel $p$ defines a linear map $\Cc^\infty(M) \to \sD'(M)$, and so these operators do not form an algebra.  In this respect, the proper distributions are better.

\begin{proposition}
\label{prop:proper_kernels}
 The following are equivalent for a distribution $p\in \sD'(M\times M)$:
 \begin{enumerate}
  \item[(1)] The Schwartz kernel operator with kernel $p$ extends to an operator preserving each of the spaces $\sD'(M)$, $\sE'(M)$, $C^\infty(M)$ and $\Cc^\infty(M)$;
  \item[(2)] $p$ is properly supported and semiregular in both variables (see Section \ref{sec:calculi});
  \item[(3)] $p\in\sDp'(M\times M)$.
 \end{enumerate}
\end{proposition}

In a similar fashion, the convolution product of distributions on a groupoid $G$,
\[
 p*q(\gamma) = \int_{\beta\in G^{r(\gamma)}} p(\beta ) q(\beta^{-1}\gamma),
\]
does not make sense for arbitrary distributions $p,q\in\sD'(G)$, but it does make sense for $p,q\in\sE'_r(G)$, via the interpretation
\[
 \ip{p*q,\varphi} = \ip{p(\beta), \ip{q(\beta^{-1}\gamma), \varphi(\gamma)}}, \qquad\qquad 
  \varphi \in C^\infty(G).
\]
(See \cite{LesManVas, VanYun:PsiDOs}.)  Moreover, this product restricts to the subspace $\sE'_{r,s}(G)$.  In particular, under the equivalences of Proposition \ref{prop:proper_kernels}, composition of Schwartz kernel operators corresponds to groupoid convolution on $\sE'_{r,s}(M\times M) \cong \sDp'(M\times M)$.  

\medskip

We next need to introduce the Lie groupoid analogue of the ideal of smoothing operators on $M$.

\begin{definition}
 Let $\Omega_r$ denote the bundle of $1$-densities along $\ker(dr)$, the tangent bundle of the $r$-fibres.  
 We denote by $\Cp^\infty(G;\Omega_r)$ the space of sections $f$ of $\Omega_r$ with \emph{proper support}, meaning that $r,s : \supp(f) \to G^{(0)}$ are both proper maps.
\end{definition}

\begin{lemma}
 The space $\Cp^\infty(G;\Omega_r)$ is a two-sided ideal of $\sE_{r,s}'(G)$.
\end{lemma}


\begin{remark}
 Note that the subalgebra of smooth densities in $\sE_r'(G)$ is a right ideal but not a left ideal.  Likewise, the algebra of smooth densities in $\sE_s'(G)$ is a left ideal but not a right ideal.  

 When $G = M \times M$, this corresponds to the fact that the properly supported smoothing operators on $M$ are a left ideal in the algebra of continuous operators on $\sE'(M)$, but not a right ideal, and a right ideal in the algebra of continuous operators on $\Cc^\infty(M)$ but not a left ideal.  For this reason, it is crucial that we work with \emph{proper} $r$-fibred distributions on $G$.
\end{remark}

\section{The osculating groupoid}

Note that the tangent bundle $TM$ of a smooth manifold is a Lie groupoid, viewed as a bundle of abelian Lie groups.
The osculating groupoid $\THM$ is a replacement for $TM$ in the world of filtered manifolds.  As usual, we start with the Lie algebroid.

\begin{definition}
 The \emph{osculating Lie algebroid} $\tHM$ of a filtered manifold $M$ is the associated graded bundle of the filtered tangent bundle:
 \[
  \tHM = \gr TM = \bigoplus_{m\in\NN} H^m / H^{m-1}.
 \]
\end{definition}

This is indeed a Lie algebroid, with anchor zero, thanks to the following calculation.  
Let $X\in \Gamma^\infty(H^m)$, $Y\in\Gamma^\infty(H^n)$ be vector fields on $M$ of $H$-order $m$ and $n$ respectively.  Then for any $f,g\in C^\infty(M)$ we have
\begin{align*}
 [fX,gY] &= fg[X,Y] + f(Xg)Y - g(Yf) X \\
         &\equiv fg[X,Y] \qquad \mod \Gamma^\infty(H^{m+n-1}).
\end{align*}
Thus the Lie bracket of vector fields induces a $C^\infty(M)$-linear Lie bracket on sections of the associated graded bundle, and hence a pointwise Lie bracket on the fibres $\tHM_x$ ($x\in M$).  In this way, $\tHM$ is a smooth bundle of nilpotent Lie algebras.  

\begin{definition}
 The \emph{osculating groupoid} $\THM$ of a filtered manifold $M$ is the bundle of connected, simply connected nilpotent Lie groups which integrates $\tHM$.
\end{definition}

In other words, $\THM_x = \tHM_x$ with product defined by the Baker-Campbell-Hausdorff formula.  

\medskip

The principal part $\sigma_m(P)$ of a differential operator of $H$-order $m$ (see Definition \ref{def:principal_part}) can now be interpreted in groupoid language.  Let us spell this out.

The classical tangent bundle $TM \onto M$ is the Lie algebroid of the pair groupoid $M \times M$.  Its section space $\Gamma^\infty(TM)$ is the space of vector fields, and its universal enveloping algebra $\sU(TM)$ is the algebra $\DO(M)$ of linear differential operators on $M$ (see \cite{NisWeiXu}).  A Lie filtration on $M$ is equivalent to a Lie algebroid filtration on $TM$, and this induces an algebra filtration on the enveloping algebra $DO(M)$.  This is precisely the filtration $\DO^\bullet_H(M)$ by $H$-order.

Passing to the associated graded spaces is functorial for all these constructions.  The osculating groupoid $\tHM = \gr TM$ is the associated graded of $TM$, and thus the universal enveloping algebra of $\tHM$ is $\sU(\tHM) = \gr\DO_H^\bullet(M)$.  The principal part of a differential operator of $H$-order $m$ is given by the canonical grading map
\begin{equation}
\label{eq:principal_part}
 \sigma_m : \DO_H^m(M) \to \sU^m(\tHM).
\end{equation}

\section{The tangent groupoid of a filtered manifold}

The tangent groupoid $\TTHM$ of a filtered manifold $M$ was introduced independently in \cite{ChoPon} and \cite{VanYun:groupoid}.  It is a smooth one-parameter family of groupoids which deforms the pair groupoid \[M\times M \rightrightarrows M\] to the osculating groupoid \[\THM\rightrightarrows M.\]  Algebraically, we have
\begin{equation}
  \label{eq:TTHM}
 \TTHM = (\THM \otimes \{0\}) \sqcup (M\times M \times \RR^\times) \quad\rightrightarrows\quad M \times \RR.
\end{equation}
The Lie groupoid structures on the two disjoint components in \eqref{eq:TTHM} are the standard ones, namely two elements of $M\times M \times \RR^\times$ are composable if and only if their $\RR^\times$ components are equal and their $M \times M$ components are composable in the pair groupoid, and two elements of $\THM \times \{0\}$ are composable if and only if their $\THM$ components are.

The difficulty lies in giving the correct global smooth structure on $\TTHM$.
As usual, this is most easily achieved by beginning with  the Lie algebroid $\ttHM$.  Note that $\ttHM$ is a deformation of $TM$ to its associated graded $\tHM = \gr(TM)$.  

\begin{notation}
 When dealing with a space $\mathfrak{X}$ which is fibred over $\RR$, we will write $\mathfrak{X}|_t$ for the fibre over $t$.  Likewise if $f$ is a function or section of a bundle over $\mathfrak{X}$ we write $f|_t$ for its restriction to $\mathfrak{X}|_t$.
\end{notation}

\begin{lemma}
 \label{lem:smooth_structure}
 There is a unique $C^\infty$ structure on the disjoint union
 \[
  \ttHM = (\tHM \times\{0\}) \sqcup (TM \times \RR^\times)
 \]
 which makes it into a smooth vector bundle over $M \times \RR$ with the the following property: for any smooth vector field $X \in \Gamma^\infty(H^m)$, the section $\XX:M\times\RR \to \ttHM$ defined by
 \begin{equation}
  \label{eq:XX}
  \XX(x,t) = \begin{cases}
              t^m X(x), & \text{if } t \neq 0, \\
              \sigma_m(X(x)), & \text{if } t= 0
             \end{cases}
 \end{equation}
 is a smooth section of $\ttHM \onto M \times \RR$.  With this smooth structure, $\ttHM$ becomes a Lie algebroid with anchor and Lie bracket defined fibrewise on each $\ttHM|_t$.
\end{lemma}

\begin{proof}
 See \cite{VanYun:groupoid}
\end{proof}

The sections $\XX$ will play a major role in the sequel.  We will write $\bbsigma_m$ for the map
\begin{equation}
 \label{eq:bbsigma}
 \bbsigma_m : \Gamma^\infty(H^m) \to \Gamma^\infty(\ttHM); \qquad X \mapsto \XX.
\end{equation}

Having defined a smooth Lie algebroid structure on $\ttHM$, there are many results which permit us to integrate it to a Lie groupoid \cite{Debord:Integration, CraFer}.  In this case, we already have smooth groupoid structures on the two components $\tHM \times \{0\}$ and $TM \times \RR^\times$, so the most easily applicable result from the literature is a theorem of Nistor \cite{Nistor:integration} (with its correction in \cite{BenNis}).  Modulo  the technical issue of $s$-simply connectedness---which we will sweep under the rug here---we can immediately deduce the existence of a $C^\infty$-structure on $\TTHM$ making it into a Lie groupoid with Lie algebroid $\ttHM$.  See \cite{VanYun:groupoid} for full details.

\section{The $\RR^\times_+$-action}

The tangent groupoid $\TTHM$ admits a crucial extra piece of structure, namely an action of $\RR^\times_+$ by Lie groupoid automorphisms, \emph{cf}.~\cite{DebSka:Rx-action}.  We begin with an $\RR^\times_+$-action on the bundle of osculating groups $\THM$, which generalizes the action of $\RR^\times_+$ on $TM$ by homotheties.

\begin{definition}
 Let $V = \bigoplus_{m\in\NN} V_m$ be a finite dimensional $\NN$-graded vector space (or vector bundle).  The \emph{dilations} $\delta_\lambda$ of $V$ are the linear automorphisms defined for each $\lambda\in\RR^\times_+$ by
 \[
  \delta_\lambda v = \lambda^m v \qquad \text{for all } v\in V_m.
 \]
\end{definition}

\begin{lemma}
 Let $M$ be a filtered manifold.  The dilations $(\delta_\lambda)_{\lambda\in\RR^\times_+}$ of $\tHM = \gr(TM)$ define a smooth action  of $\RR^\times_+$ by Lie algebroid automorphisms.  This integrates to a smooth $\RR^\times_+$-action  by groupoid automorphisms on the osculating groupoid $\THM$, which we again denote by $\delta_\lambda$.
\end{lemma}

\begin{definition}
 We define an $\RR^\times$-action $\alpha$ on $\TTHM$ by
 \begin{align*}
  \alpha_\lambda : & (x,y,t) \mapsto (x,y,\lambda^{-1}t), && t\neq 0,~ (x,y) \in M \times M \\
                   & (x,\xi,0) \mapsto (x, \delta_\lambda \xi, 0), && t=0, ~ \xi \in \THM_x.
 \end{align*}
\end{definition}

We also use the same notation $\alpha_\lambda$ for the corresponding automorphisms of the Lie algebroid $\ttHM$:
\begin{align*}
 \alpha_\lambda : & (x,v,t) \mapsto (x, v, \lambda^{-1}t), && t\neq 0,~ v\in TM_x \\
    &(x, \xi, 0) \mapsto (x, \delta_\lambda \xi , 0), && t=0,~ \xi \in \tHM_x.
\end{align*}
Note that the sections $\XX = \bbsigma(X)$ of $\ttHM$ which we defined in Lemma \ref{lem:smooth_structure} are homogeneous with respect to this action, in the sense that
\begin{equation}
 \label{eq:XX_homogeneity}
 \alpha_{\lambda*} \XX = \lambda^m \XX , \qquad \text{ for all } \lambda\in\RR^\times_+.
\end{equation}
It follows that the maps $\alpha_\lambda$ are smooth, both on the Lie algebroid $\ttHM$ and the tangent groupoid $\TTHM$.

Since the $\alpha$-action on $\TTHM$ is by Lie groupoid automorphisms, it induces an action on the spaces of $r$- and $s$-fibred distributions
\[
 \alpha_{\lambda*} : \sE'_r(\TTHM) \to \sE'_r(\TTHM), \qquad
 \alpha_{\lambda*} : \sE'_s(\TTHM) \to \sE'_s(\TTHM).
\]
These actions restrict to the subalgebras of proper $r$- and $s$-fibred distributions, and to the ideals $\Cp^\infty(\TTHM;\Omega_r)$ and $\Cp^\infty(\TTHM;\Omega_s)$ of properly supported smooth $r$- and $s$-fibred $1$-densities.


\section{Pseudodifferential operators}
\label{sec:H-PsiDOs}

We can now give the characterization of pseudodifferential operators on a filtered manifold $M$.

Let $\pp\in\sE'_r(\TTHM)$ be a properly supported $r$-fibred distribution on the tangent groupoid.
Recall that we will write $\pp|_t$ for the restriction of $\pp \in \sE'_r(\TTHM)$ to the fibre $\TTHM|_t$.  In particular,
\[
 \pp|_1 \in \sE'_r(M\times M), \qquad \pp|_0 \in \sE'_r(\THM).
\]
Thus $\pp$ is a smooth deformation of a properly supported Schwartz kernel, semiregular in both variables, to a smooth family of properly supported distributions on the osculating groups.

\begin{definition}
 Let $M$ be a filtered manifold.  An $r$-fibred distribution $\pp\in\sE'_r(\TTHM)$ is called 
 \begin{itemize}
  \item \emph{homogeneous of weight $m$} if $\alpha_{\lambda*}\pp = \lambda^m\pp$ for all $\lambda\in \RR^\times_+$;
  \item \emph{essentially homogeneous of weight $m$} if $\alpha_{\lambda*}\pp - \lambda^m\pp \in \Cp^\infty(\TTHM;\Omega_r)$ for all $\lambda\in \RR^\times_+$.
 \end{itemize}
\end{definition}

The crucial definition is the following.

\begin{definition}
 \label{def:H-PsiDOs}
 Let $p \in \sE'_r(M \times M)$ be a properly supported Schwartz kernel on $M$, semiregular in the first variable.  We say $p$ is an \emph{$H$-pseudodifferential kernel of order $\leq m$} if $p = \pp|_1$  for some properly supported $\pp\in \sE'_r(\TTHM)$ which is essentially homogeneous of order $\leq m$.  

 The associated Schwartz kernel operator $P : C^\infty(M) \to C^\infty(M)$ will be called an \emph{$H$-pseudodifferential operator}.
\end{definition}

\begin{remark}
 The proper support condition implies that $P$ restricts to an operator on $\Cc^\infty(M)$.  Moreover, we will shortly see that this definition forces $p$ to be proper, so that $P$ also extends to  $\sD'(M)$ and $\sE'(M)$ by Proposition \ref{prop:proper_kernels}.
\end{remark}

Amongst the $H$-pseudodifferential kernels are those $p=\pp|_1$ where $\pp\in\sE_r'(\TTHM)$ is homogeneous of weight $m$, and not just essentially homogeneous.  In this case, $P$ is a differential operator.

\begin{proposition}
 \label{prop:DOs}
 Let $p \in \sE'_r(M\times M)$.  Then  $p = \pp|_1$ for some $\pp \in \sE'_r(\TTHM)$ which is homogeneous of weight $m$ if and only if $p$ is the Schwartz kernel of a differential operator of $H$-order $\leq m$ on $M$.
\end{proposition}

\begin{proof}
 Let us write $\sE'_r(G)^{(0)}$ for the $r$-fibred distributions on a Lie groupoid $G$ with support on the unit space, and note that $\sE'_r(G)^{(0)}$ is isomorphic, as an algebra, to the universal enveloping algebra of the Lie algebroid  $\sA G$, see \cite{NisWeiXu, VanYun:PsiDOs}.  This isomorphism sends the homogeneous section $\XX = \bbsigma(X) \in \Gamma^\infty(\ttHM)$ from Lemma \ref{lem:smooth_structure} to a homogeneous $r$-fibred distribution in $\sE'_r(\TTHM)^{(0)}$.  This proves that the Schwartz kernels of vector fields extend to homogeneous elements of $\sE'_r(\TTHM)$.  The result extends to $\DO_H^\bullet(M)$ by multiplicativity.
 
 Conversely, note that the only orbits of $\alpha$ on $\TTHM$ which are $r$-proper are the orbits contained in the unit space.   Therefore, if $\pp$ is $r$-properly supported and homogeneous of weight $m$, we must have $\supp(\pp) \subset G^{(0)}$. The result follows.
\end{proof}

We likewise see that the singular support of an essentially homogeneous $r$-fibred distribution on $\TTHM$ is invariant for the $\RR^\times_+$-action, and so contained in $\TTHM^{(0)}$. 

\begin{proposition}
 If $\pp\in\sE'_r(\TTHM)$ is essentially homogeneous of weight $r$, then $\singsupp(\pp) \subseteq \TTHM^{(0)}$.  Therefore, $H$-pseudodifferential operators are pseudolocal.
\end{proposition}

We have thus obtained the first two of the five desirable properties (Definition \ref{def:PsiDOs}) for an algebra of pseudodifferential operators.  The other three require some more serious analysis, which we will outline in Section \ref{sec:conormality}. 

But let us pre-empt this by stating the relation of our $H$-pseudodifferential calculus with the classical calculus in the case of an unfiltered manifold.  The proof depends upon the machinery to follow.

\begin{theorem}
  Let $M$ be a manifold with the trivial filtration.  Then $\Psi^m_H(M)$ is the space of Schwartz kernels of properly supported classical (polyhomogeneous step one) pseudodifferential operators on $M$ of order $\leq m$. 
\end{theorem}

\section{Principal symbols}
\label{sec:cosymbols}

\begin{definition}
 We introduce the notation $\PPsi^m_H(M)$ for the set of $r$-fibred distributions $\pp\in\sE'_r(\TTHM)$ that are properly supported and essentially homogeneous of weight $m$.
 Thus $\Psi^m_H(M) = \PPsi^m_H(M)|_1$.  
\end{definition}

Next, we consider the restriction $\PPsi^m_H(M)|_0$.

\begin{lemma}
	\label{lem:cosymbol}
  If $\pp$ and $\pp' \in \PPsi^m_H(M)$  satisfy $\pp|_1 = \pp'|_1$, then $\pp|_0 - \pp'|_0 \in \Cp^\infty(\THM;\Omega_r)$.
\end{lemma}

\begin{proof}
 See \cite{VanYun:PsiDOs}. 
\end{proof}

The suggests the following definition.

\begin{definition}
	The space of \emph{principal cosymbols} is defined as
	\[
	 \Sigma_H^m(M) := \{ a \in \sE'_r(\THM) / \Cp^\infty(\THM;\Omega_r)
	   \st \delta_{\lambda*} a = \lambda^m a \quad \forall \lambda\in\RR^\times_+\}.
	\]
\end{definition}

To avoid clutter in the notation, we will typically use the same notation to denote a principal cosymbol class $a\in\Sigma^m_H(M)$ and an $r$-fibred distribution $a\in \sE'_r(\THM)$ which represents it.

\begin{definition}
\label{def:cosymbol_map}
Lemma \ref{lem:cosymbol} shows that we have a well-defined map
\[
  \sigma_m : \Psi^m_H(M)\to \Sigma_H^m(M); \qquad \pp|_1 \mapsto \pp|_0.
\]
This is called the \emph{principal cosymbol} map.
\end{definition}

\begin{remark}
We had previously used the notation $\sigma_m$ to denote the principal part of a differential operator,
\[
 \sigma_m : \DO^m_H(M)  = \sU^m(TM)  \to \sU^m(\tHM).
\]
Under the  isomorphism $\sU(\sA G) \cong \sE'_r(G)^{(0)}$ from the proof of Proposition \ref{prop:DOs}, this becomes a map
\[
 \sigma_m : \sE_r'(M\times M)^{(0)} \to \sE_r'(\THM)^{(0)}.
\]
In this way, one can check that Definition \ref{def:cosymbol_map} is consistent with the earlier notation.
\end{remark}

\begin{proposition}
	For every $m\in\RR$, we have a short exact sequence
	\begin{equation}
	\label{eq:ses}
	 0 \longrightarrow \Psi^{m-1}_H(M) \longrightarrow \Psi^m_H(M) \stackrel{\sigma_m}{\longrightarrow} \Sigma_H^m \longrightarrow 0.
	\end{equation}
\end{proposition}

\begin{proof}
Multiplication by the real parameter $t$ yields a linear map
\[
  \times t : \PPsi^{m-1}_H(M) \to \PPsi^m_H(M).
\]
The restriction of this map to $t=1$ yields the inclusion on the left of \eqref{eq:ses}.  For exactness in the middle, suppose $p\in\Psi^m(M)$ has $\sigma_m(p) = 0$.  Then we can find $\pp\in\PPsi^m_H(M)$ with $\pp|_1=p$ and $\pp|_0 = 0$.  Since $\pp$ is a smooth family of distributions on the $r$-fibres of $\TTHM$, $t^{-1}\pp$ is well-defined and is essentially homogeneous of weight $m-1$, so $p\in\Psi^{m-1}_H(M)$.  Finally, a construction in groupoid exponential coordinates on $\TTHM$ shows that any $a \in \Sigma^m_H(M)$ can be extended to $\mathbb{a}\in\PPsi^m_H(M)$ with $\mathbb{a}|_0 = a$, which proves the surjectivity of $\sigma_m$.
\end{proof}

\section{Algebra structure, conormality, regularity}
\label{sec:conormality}

Let $\pp, \qq \in \sE'_r(\TTHM)$.  If $\pp$ is homogeneous of weight $m$ and $\qq$ is homogeneous of weight $n$, then a simple calculation shows that their convolution product $\pp*\qq$ is homogeneous of weight $m+n$.  It is tempting to say the same of \emph{essentially} homogeneous elements, but we are obstructed by the fact that $\Cp^\infty(\TTHM;\Omega_r)$ is only a right ideal in $\sE'_r(\TTHM)$, not a two-sided ideal.  

For this reason, the following technical lemma is crucially important for further progress.  Its proof, which also serves to give the regularity results to follow, is one of main pieces of analysis in this work.  We won't give complete details, but we will sketch the main ideas.  For full details, see \cite{VanYun:PsiDOs}.

\begin{lemma}
	\label{lem:properly fibred}
        Let $\pp\in\sE'_r(\TTHM)$ be properly supported.  If $\pp$ is essentially homogeneous of weight $m$ for some $m$, then $\pp\in \sE'_{r,s}(\TTHM)$.
\end{lemma}

\begin{proof}[Sketch of proof]
	\begin{itemize}
		\item 	We first linearize the problem.  The groupoid exponential allows us to identify a neighbourhood of the unit space $\TTHM^{(0)}$ with a neighbourhood of the zero section in $\ttHM$.  By fixing a splitting $\tHM = \gr(TM) \stackrel{\cong}{\to} TM$, this can in turn be identified with a neighbourhood of the zero section in $\tHM \times \RR$.  Under this identification, $\pp$ is identified with a smooth family of compactly supported distributions on the fibres of the vector bundle $\tHM \times \RR \onto M\times \RR$
                which is essentially homogeneous with respect to the $\RR^\times_+$-action
		\[
		 \alpha_\lambda:\tHM\times\RR \to \tHM\times\RR; \qquad
		   (x,\xi,t) \mapsto (x,\delta_\lambda\xi,\lambda^{-1}t).
		\]

		\item  Now apply the fibrewise Fourier transform:
		\[
		 \widehat~ : \sE'_r(\tHM\times\RR) \to C^\infty(\tsHM\times\RR); \qquad \pp \mapsto \hat\pp,
		\]
                where $\tsHM \onto M $ is the dual bundle of $\tHM$.
		Then $\hat\pp$ is a smooth function on $\tsHM\times \RR$ with its own essential homogeneity property.  Specifically, let $\delta'_\lambda$ denote the canonical dilations on $\tsHM$ with $\NN$-grading dual to that of $\tHM$, and let $\beta_\lambda$ denote the action on $\tsHM\times \RR$ given by
		\[
		\beta_\lambda:\tsHM\times\RR \to \tsHM\times\RR; \qquad
		(x,\eta,t) \mapsto (x,\delta'_\lambda\eta,\lambda t).
		\]
    Then for all $\lambda\in\RR^\times_+$ we have
		\[
		  \beta_\lambda^* \hat\pp - \lambda^m \hat\pp \in \sS_r(\tsHM \times \RR),
		\]
		where $\sS_r(\tsHM \times \RR)$ denotes the smooth functions on $\tsHM\times\RR$ which are Schwartz-class on each $r$-fibre (see \cite{VanYun:PsiDOs} for the precise definition).
		
		\item Essential homogeneity of smooth functions near infinity is better behaved than essential homogeneity of distributions near $0$.  Specifically, a bundle version of \cite[Lemma 2.2]{Taylor:microlocal} shows that, outside any neighbourhood of the zero section, $\hat\pp$ is equal to a genuinely homogeneous smooth function plus a function of rapid decay on the fibres.  As a consequence, we obtain bounds on all derivatives of $\hat{\pp}(x,\eta,t)$ in terms of powers of $\|\eta\|$.  
		
		\item The above bounds on $\hat\pp$ imply that the wavefront set of $\pp$ is conormal to the unit space.  This is sufficient to deduce that $\pp$ is both $r$- and $s$-fibred, thanks to \cite[Proposition 7]{LesManVas}.
	\end{itemize}
\end{proof}

\begin{corollary}
 Convolution of $r$-fibred distributions induces an $\RR$-filtered algebra structure on $\Psi^\bullet_H(M)$.
\end{corollary}

\begin{proof}
	We have $\PPsi^m_H(M) \subseteq \sE'_{r,s}(\TTHM)$.  Since $\Cp^\infty(\TTHM;\Omega_r)$ is a two-sided ideal in $\sE'_{r,s}(\TTHM)$, the convolution product induces maps $\PPsi^m_H(M) \times \PPsi^n_H(M) \to \PPsi^{m+n}_H(M)$ for every $m,n\in\RR$.  Restriction to $t=1$ gives the product on $\Psi^\bullet_H(M)$.
\end{proof}

Going further, a careful analysis of the bounds on $\hat\pp$ allows us to deduce the degree of regularity of $\pp$ on the unit space of $\TTHM$.

\begin{lemma}
 \label{lem:regularity}
	Let $N$ be the depth of the filtration on $TM$, and let $d_H = \linebreak \sum_{k=1}^N k\dim(H^k/H^{k-1})$ be the \emph{homogeneous dimension} of $\tHM$.  Then
	\[
	 \Psi^{-d_H-kN-1}_H(M) \subseteq \Cp^k(M\times M; \Omega_r).
	\]
\end{lemma}

\begin{proof}
 See \cite{VanYun:PsiDOs}.
\end{proof}

\begin{corollary}
	We have $\bigcap_{m\in\RR} \Psi^m_H(M) = \Psi^\infty(M)$, the algebra of smooth kernels.
\end{corollary}

This is Property (3) from our list in Definition \ref{def:PsiDOs}.  We also obtain Property (4), the existence of asymptotic limits.

\begin{corollary}
	\label{cor:asmptotic_limits}
	Given any sequence $(p_i)_{i\in\NN}$ of $H$-pseudodifferential kernels, with $p_i \in \Psi^{m-i}_H(M)$ for each $i\in\NN$, there exists $p\in\Psi^m_H(M)$ such that for any $l\in\NN$,
	\[
	  p - \sum_{i=0}^l p_i \in \Psi^{m-l-1}_H(M).
	\]
\end{corollary}

\begin{proof}
 Let $\pp_i\in \PPsi^{m-i}_H(M)$ with $\pp_i|_1 = p_i$.   Lemma \ref{lem:regularity} implies that for $i$ sufficiently large we have $\pp_i \in \Cp^{k(i)}(\TTHM)$ for some $k(i)$ with $k(i) \to \infty$ as $i\to \infty$.  After modifying each $\pp_i$ by an element of $\Cp^\infty(\TTHM;\Omega_r)$, we may assume that $\pp_i$ vanishes on the unit space $\TTHM^{(0)}$ to order $k(i)$.  
 
 Then we have $t^i\pp_i \in \PPsi^m_H(M)$, vanishing to order $i$ on $\TTHM|_0$ and to order $k(i)$ on $\TTHM^{(0)}$.  An analogue of Borel's Lemma yields $\pp\in \sE'_r(\TTHM)$ such that $\pp - \sum_{i=0}^l t^{i}\pp_i \in t^{l+1}\Cp^{k(l+1)}$ for all $l\gg 0$.  Since the finite sum $\sum_{i=0}^l t^{i}\pp_i$ is homogeneous of weight $m$ modulo $\Cp^\infty(\TTHM;\Omega_r)$, it follows that $\pp$ is homogeneous of weight $m$ modulo $\Cp^k(\TTHM ;\Omega_r)$ for every $k$, and hence modulo $\Cp^\infty(\TTHM ;\Omega_r)$.  That is, $\pp \in \PPsi^m_H(M)$.    Put $p=\pp|_1$.  Then $t^{-l-1}(\pp - \sum_{i=0}^l t^{i}\pp)\in \PPsi^{m-l-1}_H(M)$, so $p - \sum_{i=1}^l p_i \in \Psi_H^{m-l-1}$.
 This completes the proof. 
\end{proof}

\section{Hypoellipticity}
\label{sec:Rockland}

Given the above structure, the following definition is the natural abstract analogue of ellipticity for $H$-pseudodifferential operators.

\begin{definition}
	An $H$-pseudodifferential kernel $p\in\Psi^m_H(M)$ is called \emph{$H$-elliptic} if its principal cosymbol $\sigma_m(p)$ admits a convolution inverse in $\Sigma^{-m}_H(M)$.
\end{definition}

By the surjectivity of the principal cosymbol map, this means that there exists $q\in\Psi^{-m}$ such that 
\[
p*q -I \in \Psi^{-1}_H(M).
\]
Thus, $H$-elliptic operators satisfy axiom (5) of Definition \ref{def:PsiDOs}. 
The following is now an immediate consequence of Theorem \ref{thm:Ell}.

\begin{theorem}
	Any $H$-elliptic operators on a filtered manifold $M$ is hypoelliptic. 
\end{theorem}

The $H$-ellipticity condition is difficult to verify in practice.  The correct analogue for ellipticity on a filtered manifold should be the Rockland property \cite{Rockland, HelNou, Melin:preprint}.  The relationship between these two properties has been clarified by Dave and Haller \cite{DavHal:BGG}, based on results of Folland-Stein \cite{FolSte:Hardy} and Ponge \cite{Ponge:Heisenberg}.  We begin by recalling the definition.

\begin{definition}
	\label{def:Rockland}
	A differential operator $P\in\DO_H^m(M)$ on a filtered manifold $M$ is \emph{Rockland} if, for every $x\in M$ and for every non-trivial irreducible unitary representation $\pi$ of the osculating group $\tHM_x$ on a Hilbert space $H^\pi$, the operator $\pi(\sigma_m(P))$ is injective on the space of smooth vectors in $H^\pi$.  We will say \emph{$P$ is two-sided Rockland} if both $P$ and its transpose $P^t$ are Rockland.
\end{definition}

\begin{theorem} \textup{\cite{DavHal:BGG}}
	If $P\in \DO(M)$ two-sided Rockland then it is $H$-elliptic.  Therefore, every Rockland operator on a filtered manifold is hypoelliptic.
\end{theorem}

It is also shown in \cite{DavHal:BGG} that there is a Sobolev theory adapted to the pseudodifferential calculus $\Psi_H^\bullet$.   It follows that if $M$ is compact without boundary, then $H$-elliptic operators are Fredholm operators between the appropriate Sobolev spaces \cite[Corollary 3.28]{DavHal:BGG}.

\chapter{Pseudodifferential operators on multifiltered manifolds and equivariant index theory for complex semisimple groups}
\label{sec:flags}

Flag manifolds of semisimple Lie groups $G$ are important examples of filtered manifolds.  In differential geometry, they are the model spaces of a class of spaces called \emph{parabolic geometries}, see \emph{e.g.}\ \cite{CapSlo}, which generalize CR manifolds and play an important role in twistor theory \cite{BasEas}, amongst others.  Our interest, however, will come from equivariant index theory.

Kasparov and Julg \cite{JulKas, Julg:Spn1} showed that when $G=\SU(n,1)$ or $\Sp(n,1)$, one can use natural subelliptic operators on the flag manifold of $G$ to define a $G$-equivariant $K$-homology class which represents the notorious element $\gamma \in \KK^G(\CC,\CC)$ at the heart of the Baum-Connes conjecture.  On the other hand, when $G$ has real rank greater than $1$, a result of Puschnigg \cite{Puschnigg:T} implies that any construction based on the filtered calculus of Chapter \ref{sec:filtered_manifolds} can only result in multiples of the trivial class. 

Nonetheless, at least in the case of $G=\SL(3,\CC)$, it is possible to represent the $\gamma$ element as an equivariant $K$-homology class for the flag manifold by using the Bernstein-Gelfand-Gelfand complex \cite{Yuncken:products_of_PsiDOs, Yuncken:BGG}.  The Bernstein-Gelfand-Gelfand (BGG) complex is a $G$-equivariant differential complex built from longitudinally elliptic differential operators along the fibres of the numerous canonical fibrations of the flag manifold.
The analysis in those articles was based on Gelfand-Tsetlin theory for the representations of $\SU(3)$, but these methods generalize badly, since analogues of Gelfand-Tsetlin theory for Lie algebras other than type $A_n$ are difficult.

Here we shall focus on the article \cite{Yuncken:foliations}, which describes a rudimentary pseudodifferential theory adapted to the analysis of longitudinal differential operators on manifolds with multiple foliations.  It is rudimentary in the sense that it essentially only distinguishes between operators of order $0$ and operators of negative order.  Nevertheless, this suffices to distinguish bounded and compact operators, and their analogues for foliations, which will be enough for application to index theory.


\section{Complex flag manifolds}
\label{sec:flag_notation}

Throughout, $G$ will denote a connected, simply connected complex semisimple Lie group.  We will use the following notation:
\begin{itemize}
 \item $\lie{g}$ is the Lie algebra of $G$
 \item $\lie{h}$ is a Cartan subalgebra
 \item $\bDelta$, $\bDelta^+$, $\bSigma$ denote the sets of roots, positive roots, and simple roots, respectively.
 \item $\rho = \half\sum_{\beta\in\bDelta^+}\beta$ is the half-sum of the positive roots,
 \item $\bP$ and $\bQ$ are the integral weight and root lattices, respectively.
 \item $\lie{k}$, $\lie{a}$, $\lie{n}$ are the components of the Iwasawa decomposition, with $\lie{n} = \bigoplus_{\alpha\in\bDelta^+}\lie{g}_\alpha$ the sum of the positive root spaces.
 \item $\bar{\lie{n}}= \bigoplus_{\alpha\in\bDelta^+}\lie{g}_{-\alpha}$ is the opposite nilpotent subgroup.
 \item $\lie{b} = \lie{h}\oplus \lie{n}$ and $\bar{\lie{b}} = \lie{h}\oplus \bar{\lie{n}}$ are the standard upper and lower Borel subgroups.
\end{itemize}
A connected Lie subgroup of $G$ will always be denoted by the same letter as its Lie algebra, in upper case.

The standard parabolic Lie subalgebras of $\lie{g}$ are in bijection with subsets $I\subseteq \bSigma$ of the simple roots. To make this precise this, we introduce the notation
\[
 \ip{I} = \{\beta \in \bDelta^+ \st \beta = \sum_{\alpha \in I} n_\alpha \alpha \text{ for some } n_\alpha \geq 0\}
\]
for the set of positive roots in the span of $I$.  Then we put%
\footnote{
Unfortunately, our notation here is the opposite of standard one, in that we are writing $\lie{p}_I$ for what would often be denoted $\lie{p}_{\bSigma\setminus I}$.  We choose this convention to be consistent with \cite{Yuncken:BGG, Yuncken:foliations}.
}
\[
 \lie{p}_I = \lie{b} \oplus \sum_{\alpha \in \ip{I}} \lie{g}_{-\alpha}.
\]
In particular, $\lie{p}_\bSigma = \lie{g}$ and $\lie{p}_\emptyset = \lie{b}$.  The homogeneous spaces
\[
 \sX_I = G / P_I
\]
are the \emph{flag manifolds} of $G$.  

In the case where $I= \{\alpha\}$ is a singleton, we will write $\lie{p}_\alpha$, $P_\alpha$, $\sX_\alpha$, \emph{etc.}\ for $\lie{p}_{\{\alpha\}}$, $P_{\{\alpha\}}$, $\sX_{\{\alpha\}}$, \emph{etc}.
We will also write $\sX = \sX_\emptyset = G/B$ for the full flag manifold.

The various flag manifolds are connected by smooth $G$-equivariant fibrations $\sX_J \onto \sX_I$ whenever $J\subseteq I$.  In particular we have $G$-equivariant maps
\[
 q_I : \sX \onto \sX_I.
\]
Thus, for each $I\subseteq \bSigma$ there is a foliation of the full flag manifold by the fibres of $q_I$.  We will denote the tangent bundle to the fibres by
\[
 \sF_I = \ker dq_I \subseteq T\sX.
\]

\section{Motivation: The BGG complex}
\label{sec:BGG}

As mentioned above, this chapter will be dedicated to describing the results of \cite{Yuncken:foliations}, concerning longitudinal pseudodifferential operators on manifolds with multiple foliations.  But this work was entirely motivated by the equivariant index theory of the differential Bernstein-Gelfand-Gelfand complex, which consists of longitudinal differential operators along the various fibrations of the flag manifold just described.  For context, therefore, we will start with a very rapid overview of the Bernstein-Gelfand-Gelfand complex.

\bigskip


Let $\mu\in\bP$ be an integral weight of $\lie{g}$, and let $\lambda \in \lie{a}_\CC^* = \lie{h}^*$.  These exponentiate to a character $e^\mu$ of the compact torus $M = K\cap H$ and a character $e^\lambda$ of $A$ (not generally unitary), via the formulas
\begin{align*}
 e^\mu(\exp(X)) &:= e^{\mu(X)}, & X\in\lie{m}, \\
 e^\mu(\exp(Y)) &:= e^{\mu(Y)}, & Y\in\lie{a}. 
\end{align*}
We obtain a character of the Borel subgroup $B=MAN$ by
\[
 \chi_{\mu,\lambda}(man) = e^\mu(m) e^\lambda(a).
\]

\begin{definition}
 \label{def:induced_bundle}
 We write
 \[
  \sE_{\mu,\lambda} = G \times_{B} \CC_{\mu,\lambda+2\rho}
 \]
 for the $G$-equivariant line bundle over the flag manifold $G/\bar{B}$ which is induced from the shifted character $\chi_{\mu,\lambda+2\rho}$.
 This means that smooth sections of $\sE_{\mu,\lambda}$ are defined by
\begin{equation}
\label{eq:induced_sections}
 \Gamma^\infty(\sE_{\mu,\lambda}) = \{ \xi \in C^\infty(G) \st
    \xi(gb) = \chi_{\mu,\lambda+2\rho}(b^{-1}) \xi(g) \quad \forall g\in G,~ b \in B  \}.
\end{equation}
\end{definition}

The space of $L^2$-sections of $\sE_{\mu,\lambda}$, which we denote by $L^2(\sE_{\mu,\lambda})$, is the completion of $\Gamma^\infty(\sE_{\mu,\lambda})$ with respect to the $K$-invariant inner product
\[
 \ip{\xi, \eta} = \int_K \overline{\xi(k)} \eta(k)\, dk.
\]
The left regular representation of $G$ on $L^2(G)$ restricts to a representation on $L^2(\sE_{\mu,\lambda})$  called the \emph{principal series representation} of parameter $(\mu,\lambda) \in \bP \times \lie{h}^*$.  It is unitary when $\lambda \in i\lie{a}^*$, see \emph{e.g.} \cite{Knapp:representation_theory}.

The \emph{Bernstein-Gelfand-Gelfand complex} is a complex of $G$-equivariant differential operators between the section spaces of certain very particular line bundles $\sE_{\mu,\lambda}$ over $\sX$.  These operators are obtained by the right action of particular elements of $\sU(\lie{g})$.  Such operators, if well-defined, obviously commute with the principal series actions of $G$.  The details are as follows.

Let $W$ denote the Weyl group of $G$.  We write $\ell(w)$ for the length of $w\in W$, \emph{i.e.}, the minimal length as a word in the simple reflections.  We recall that the Bruhat graph is a directed graph with vertices $W$ and edges $w\to w'$ whenever $\ell(w') = \ell(w)+1$ and $w'=s_\alpha w$ for some simple reflection $s_\alpha$; see, \emph{e.g.}, \cite{Humphreys:O}.  

Note that, with the parametrization of Definition \ref{def:induced_bundle}, the trivial bundle is $\sE_{0,-2\rho}$.  We have an inclusion $\CC \into \Gamma^\infty(\sX_{0,-2\rho})$ of the trivial $G$-module as constant functions.  The Bernstein-Gelfand-Gelfand complex is a resolution of this inclusion by direct sums of principal series representations.

\begin{theorem}
 \label{thm:BGG}
 Let $w,w'\in W$ be connected by an edge $w \to w'$ in the Bruhat graph.  There exists $X\in\sU(\lie{g})$ such that the right regular action of $X$ defines a $G$-equivariant differential operator
 \begin{equation}
  \label{eq:G-DO}
  X : \Gamma^\infty(\sE_{-w\rho+\rho, -w\rho-\rho}) \to \Gamma^\infty(\sE_{-w'\rho+\rho,-w'\rho-\rho}).
 \end{equation}
 In particular, if $w'=s_\alpha w$ where $s_\alpha$ is the simple reflection associated to $\alpha\in\bSigma$, then $X=E_\alpha^n$ where $E_\alpha\in\lie{g}$ is the simple root vector of weight $\alpha$ and $n\in\NN$ is determined by $w\rho - w'\rho = n\alpha$.

 With an appropriate choice of signs, these operators $X$ form a resolution of the trivial $G$-module
 \begin{equation}
  \label{eq:BGG}
 \CC \hookrightarrow
   \Gamma^\infty(\sE_{0,-2\rho}) \to \!\!
  \bigoplus_{\ell(w) = 1} \!\!\!\! \Gamma^\infty(\sE_{-w\rho+\rho, -w\rho-\rho}) \to \!\!
  \bigoplus_{\ell(w) = 2} \!\!\!\! \Gamma^\infty(\sE_{-w\rho+\rho, -w\rho-\rho}) \to \cdots 
 \end{equation}
\end{theorem}

\begin{remark}
 With other choices of parameters $(\mu,\lambda)\in\bP\times\lie{h}^*$, one can obtain more general Bernstein-Gelfand-Gelfand complexes giving resolutions of other finite dimensional $G$-modules, see \emph{e.g.} \cite{BasEas}.
\end{remark}

The key observation for what follows is that the differential operators $E_\alpha^n$ corresponding to simple edges in the Bruhat graph are longitudinally elliptic differential operators along the fibres of the fibration $\sX \onto \sX_\alpha$.  Therefore, in order to convert the BGG complex into an equivariant $KK$-cycle, we need some  kind of pseudodifferential theory which simultaneously contains longitudinal pseudodifferential operators along each of these canonical fibrations.  This is the goal of \cite{Yuncken:foliations}.

\begin{remark}
The Dolbeault complex of $\sX$ is also a $G$-equivariant resolution of the trivial module, but being elliptic it cannot easily be converted into a nontrivial $KK$-cycle, because of the previously mentioned result of Puschnigg \cite{Puschnigg:T}.  The main problem is that the action of $G$ on the spinor bundle of $\sX$ is not isometric, nor even conformal.  For the importance of conformality in equivariant $KK$-theory, compare \cite{Kasparov:Lorentz, JulKas, Julg:Spn1}.

The advantage of the BGG complex is that the induced bundles are direct sums of line bundles, each of which admits a Hermitian metric which is conformal for the $G$-action and with explicit formulas for the conformality constants.  This is a key point in \cite{Yuncken:BGG}.
\end{remark}

\section{Longitudinal pseudodifferential operators}


\begin{definition}
Let $\sE, \sE'$ be vector bundles over the full flag variety $\sX$ of $G$.  For each $I\subseteq\bSigma$, we write $\Psi^m_I(\sE,\sE')$ for the set of longitudinal pseudodifferential operators $T:\Gamma^\infty(\sE) \to \Gamma^\infty(\sE')$ which are of order $\leq m$ (polyhomogeneous step $1$) along the leaves of the foliation $\sF_I$.
\end{definition}

If $\sE=\sE'$ we will simplify the notation to $\Psi^{-m}_I(\sE)$.   But in general, we will abbreviate all of these by writing $\Psi^m_I$, where the bundles involved can be inferred from the context.  In \cite{Yuncken:foliations} we used the notation $\Psi^m(\sF_I)$ instead of $\Psi^m_I$.

For each $I$, the elements of $\Psi^0_I(\sX)$ extend to bounded operators on $L^2(\sX)$, and thus the norm closure $\overline{\Psi^0_I(\sX)}$ is a $C^*$-algebra.  The norm closure $\overline{\Psi^{-\infty}_I(\sX)}$ of the longitudinal smoothing operators is an ideal in $\overline{\Psi^0_I(\sX)}$ which contains $\Psi^m_I(\sX)$ for all $m<0$.  It is sometimes denoted
$C^*(\sF_I)$.  

In particular, if $I=\bSigma$ then the foliation of $\sX$ by $\sF_\bSigma$ consists of a single leaf, and we have $\Psi^m_\sigma(\sX) = \Psi^m(\sX)$.  Therefore $\overline{\Psi^{-\infty}(\sF_\bSigma)}=\KK(L^2(\sX))$.

If $I \neq J$, the product of an element of $\Psi^m_I(\sX)$ and one of $\Psi_J(\sX)$ is typically not a pseudodifferential operator in any reasonable sense.  In this respect, the $C^*$-closures are much better behaved.

\begin{theorem}[\cite{Yuncken:foliations}]
	\label{thm:foliations}
	For any $I,J \subseteq \bSigma$, we have 
	\[
	 \overline{\Psi^{-\infty}_I}.\overline{\Psi^{-\infty}_J}
	   \subseteq \overline{\Psi_{I\cup J}^{-\infty}}.
	\]
	In particular, if $I_1,\ldots,I_n \subseteq \bSigma$ with $\bigcup_{k=1}^n I_k=\bSigma$, and if $P_k$ is a longitudinal pseudodifferential operator of negative order along the leaves of $\sF_{I_k}$ for each $k=1,\ldots n$, then the product $P_1\dots P_n$ is a compact operator.
\end{theorem}

This theorem is a key point for the construction of an equivariant $K$-homology class from the Bernstein-Gelfand-Gelfand complex, as we shall show in Section \ref{sec:Kaparov_product}.  But first we will describe the proof of Theorem \ref{thm:foliations} in the next three sections.

\section{Multifiltered manifolds}

Let us fix some notation.  Fix $r\in\NN^\times$, usually signifying the rank of $G$.  We denote multi-indices by $\bfa = (a_1, \ldots, a_r) \in \NN^r$.  We write $\bfa\vee\bfb$ and $\bfa\wedge\bfb$ for the entry-wise maximum and entry-wise minimum, respectively, of two multi-indices.

\begin{definition}
 Let $V$ be a finite dimensional vector space (or a finite dimensional vector bundle).  An \emph{$r$-multifiltration} on $V$ will mean a family of subspaces (or subbundles) $V^\bfa \leq V$ indexed by $\bfa \in \NN^r$ satisfying
 \begin{itemize}
  \item $V^\mathbf{0} = 0$ and $V^\mathbf{m} = V$ for some $\mathbf{m}\in\NN^r$,
  \item $V^\bfa \cap V^\bfb = V^{\bfa \wedge \bfb}$ for all $\bfa,\bfb \in \NN^r$.
 \end{itemize}
\end{definition}

\begin{definition}
 \label{def:multifiltered_manifold}
 An \emph{$r$-multifiltered manifold} is a manifold $M$ equipped with a $r$-multifiltration of its tangent bundle by subbundles $H^\bfa$ ($\bfa\in\NN^r$)  such that
 \[
  [\Gamma^\infty(H^\bfa), \Gamma^\infty(H^\bfb)] \subseteq \Gamma^\infty(H^{\bfa+\bfb}), \qquad\qquad 
   \forall \bfa,\bfb\in\NN^r.
 \]
\end{definition}

\begin{example}
 The Cartesian product of two manifolds, $M = M_1 \times M_2$ is a $2$-multifiltered manifold, with
 \[
  H^{(1,0)} = \pr_1^* TM_1, \qquad H^{(0,1)} = \pr_2^* TM_2, \qquad H^{(1,1)} = TM,
 \]
 where $\pr_i : M \onto M_i$ is the coordinate projection.
 This multifiltered structure is relevant to the study of bisingular operators \cite{Rodino:bisingular, Bohlen:bisingular} and the exterior product in analytic $K$-homology.
\end{example}

\begin{example}
 \label{eq:flag_multifiltration}
 The full flag manifold $\sX$ of a complex semisimple Lie group $G$ is an $r$-multifiltered manifold, where $r$ is the rank of $G$.
 To see this, recall that the tangent space $T\sX$ is isomorphic to the induced bundle $G \times_{B} \lie{g}/\lie{b}$, where $\lie{g}/\lie{b}$ is equipped with the adjoint action of $B$.  A $G$-equivariant multifiltration on $T\sX$ can therefore be determined by a multifiltration on the ${B}$-module $\lie{g}/\lie{b}$.  

 Let $\alpha_1, \ldots, \alpha_r$ be an enumeration of the simple roots, and for each $\bfa= (a_i,\ldots,a_r)\in\NN^r$ put
 \[
  V^{\bfa} = \left( \bigoplus  \left. \left\{ \lie{g}_\beta \st \beta = \sum_{i=1}^r b_i\alpha_i \text{ with } b_i \geq -a_i ~ \forall i \right\} \right) \right/ \lie{b}
 \]
 and
 \[
  H^\bfa = G\times_{B}\! V^\bfa.
 \]
 This makes $\sX$ into an $r$-multifiltered manifold.

This multifiltration on $\sX$ is related to the tangent bundles $\sF_I$ of the canonical fibrations $\sX\onto \sX_I$ as follows.  Multi-indices $\bfa$ are in bijection with the positive elements of the root lattice, by associating $\bfa$ to $\sum_ia_i\alpha_i$.  Then $\sF_I = H^\bfb$, where $\bfb$ corresponds to the maximal $\beta\in\bDelta^+$ with support in $I$.
\end{example}

\section{Locally homogeneous structures}
\label{sec:locally_homogeneous}

In the proof of Theorem \ref{thm:foliations} we will take advantage of the fact that multifiltration of the flag manifold is locally diffeomorphic to a multifiltration of $\bar{N}$ induced by a family of nilpotent subgroups $\bar{N}_I$.  Let us begin with some abstract definitions.

We use the notation $\epsilon_i = (0,\ldots,1,\ldots,0)$ for the multi-index whose only nonzero coefficient is the $i$th, which is $1$.
\begin{definition}
 An \emph{$r$-multigraded nilpotent Lie algebra} will mean a finite dimensional nilpotent Lie algebra $\lie{n}$ with an $\NN^r$-grading
 \[
  \lie{n} = \bigoplus_{\bfa\in\NN^r} \lie{n}_\bfa
 \]
 that is compatible with the Lie bracket, such that $\lie{n}_\mathbf{0} = 0$ and $\lie{n}$ is generated as a Lie algebra by the subspaces $\lie{n}_{\epsilon_i}$.
\end{definition}

The $r$-multigrading on $\lie{n}_\bullet$ induces an $r$-multifiltration, denoted $\lie{n}^\bullet$, by putting
\[
 \lie{n}^\bfa  = \bigoplus_{\bfb \leq \bfa} \lie{n}_\bfb.
\]
Let $N$ be the associated simply connected nilpotent Lie group and identify $TN = N \times \lie{n}$ via left translations.  Then $TN$ inherits an $r$-multifiltration, and so $N$ is a multifiltered manifold.

\begin{definition}
 \label{def:locally_homogeneous}
 Let $\lie{n}$ be a multigraded Lie algebra.  A multifiltered manifold $M$ will be called \emph{locally homogeneous of type $\lie{n}$} if there exists an atlas of local charts $\varphi: U \to M$ from open subsets of $U \subseteq N$ satisfying
 \[
  \varphi^* (H^\bfa) = U \times \lie{n}^\bfa , \qquad \qquad \forall \bfa\in\NN^r,
 \]
 where the $H^\bfa$ are the subbundles defining the multifiltration of $M$, as in Definition \ref{def:multifiltered_manifold}.
\end{definition}

Of course, we are only really interested in this structure because it encapsulates the structure of flag manifolds.

\begin{example}
 \label{ex:flag_multifiltration}
 Let $\bar{\lie{n}}$ be the nilpotent radical of the lower Borel subgroup of $\lie{g}$.  Then $\bar{\lie{n}}$ is an $r$-multigraded lie algebra via
 \[
  \bar{\lie{n}}^\bfa = \lie{g}_{-\alpha},
 \]
 where, $\alpha = \sum_i a_i\alpha_i$. 

 We claim that the multifiltration of $\sX$  in Example \ref{ex:flag_multifiltration} is locally homogeneous of type $\bar{\lie{n}}$.
 For each $x\in G$, we can define a chart
 \begin{align*}
  \varphi_x :  \bar{N} \to \sX ; \qquad
   \varphi_x(\bar{n}) = x\bar{n}B,
 \end{align*}
 which is a diffeomorphism of $\bar{N}$ onto a left translate of the big Bruhat cell in $\sX$.  From Example \ref{ex:flag_multifiltration}, we see that $\varphi_x^* (H^\bfa) = \bar{N} \times \bar{\lie{n}}^\bfa$, as desired.
\end{example}

The following definition, which we will not make much use of, allows one to generalize the foliations $\sF_I$ of the flag manifold to any locally homogeneous multifiltered manifold.

\begin{definition}
	Let $\lie{n}$ be an $r$-multigraded Lie algebra.  For any set $I\subseteq\{1,\ldots,r\}$, we will denote by $\lie{n}_I$ the Lie subalgebra of $\lie{n}$ generated by $\bigoplus_{i\in I}\lie{n}_{\epsilon_i}$.  The associated Lie subgroup is denoted $N_I$.
	
	If $M$ is any locally homogeneous multifiltered manifold of type $\lie{n}$, then the foliation of $N$ by cosets of $N_I$ induces, via the local charts, a foliation of $M$.  We denote the tangent bundle to the leaves by $\sF_I$.
\end{definition}

\section{Products of longitudinally smoothing operators}

The point of the above abstract nonsense is merely to show that Theorem \ref{thm:foliations} can be reduced to a theorem about smoothing operators along the coset foliations of a family of subgroups of a nilpotent group.  Indeed, by using a partition of unity subordinate to the charts $\varphi$ of Definition  \ref{def:locally_homogeneous}, and by putting $H_1 = N_I$, $H_2 = N_J$ and $H=N_{I\cup J}$ it suffices to prove the following result.

\begin{theorem}
	\label{thm:nilpotent_foliations}
	Let $H$ be a nilpotent Lie group and let $H_1, H_2$ be connected subgroups of $H$ such that their Lie algebras $\lie{h}_1$, $\lie{h}_2$ generate the Lie algebra $\lie{h}$ of $H$.  Let $A_1, A_2:\Cc^\infty(H) \to \Cc^\infty(H)$ be compactly supported longitudinally smoothing operators along the coset foliations of $H$ by $H_1$ and $H_2$, respectively.  Then $A_1A_2$ extends to a compact operator on $L^2(H)$.
\end{theorem}

\begin{proof}[Sketch proof]
	The longitudinally smoothing operators $A_i$ are given by a convolution formula of the following form:
	\begin{equation}
		\label{eq:longitudinally_smoothing}
	 A_i u(x) = \int_{h\in H_i} a_i(x,h)u(h^{-1}x)\, dh, \qquad u\in\Cc^\infty(H),
	\end{equation}
	for some $a_i \in \Cc^\infty(H\times H_i)$.
	
	To prove that $A_1A_2$ is compact, it suffices to prove that $A_2^*A_1^*A_1A_2$ is compact, or indeed that $(A_2^*A_1^*A_1A_2)^n$ is compact for some $n\in\NN$.  Let us take $n$ large enough that the space of products
	\[
	 \underbrace{H_2H_1H_2H_1H_2\ldots H_1H_2}_{2n+1 \text{ terms}}
	\]
  has positive measure in $H$.
  Repeated application of Equation \eqref{eq:longitudinally_smoothing} shows that
  $(A_2^*A_1^*A_1A_2)^n$ has the form
  \begin{multline*}
   (A_2^*A_1^*A_1A_2)^n u(x) \\
     =\int_{\underbrace{H_2\times H_1 \times \cdots \times H_2}_{2n+1 \text{ terms}}}
     	a(x,h_1,h_2,\ldots,h_{2n+1}) u((h_1h_2\cdots h_{2n+1})^{-1}x) \\[-2ex]
     	dh_1 \, dh_2  \ldots dh_{2n+1}
  \end{multline*}
  for some smooth function $a \in \Cc^\infty(H \times \overbrace{H_2\times H_1 \times \cdots \times H_2}^{2n+1 \text{ terms}})$.
  
  The product map 
  $
    \phi: H_2 \times H_1 \times \cdots \times H_1 \to H
  $
  is real analytic, so its derivative is surjective outside a set of measure zero.  The proof is then completed by the following Lemma, which is essentially a consequence of the implicit function theorem.

\begin{lemma} \textup{(\cite{Yuncken:foliations})}
	Let $M$ be a smooth manifold, $\phi:M\to H$ a smooth map whose derivative is surjective outside a set of zero measure, and let $a\in\Cc^\infty(M\times H)$.  Then the operator $A$ defined by
	\[
	  Au(x) = \int_M a(x,m) u(\phi(m)^{-1}x) \,dx
	\]
	is a compact operator on $L^2(H)$.
\end{lemma}

\vspace{-3.5ex}
\end{proof}

\section{Application: The $KK$-class of the BGG complex of rank-two complex semisimple groups}
\label{sec:Kaparov_product}

To conclude this chapter, we return to the motivating application of constructing a $G$-equivariant $K$-homology class for the full flag manifold from the Bernstein-Gelfand-Gelfand complex.
At present, this construction has only been completely carried out for the group $G=\SL(3,\CC)$ (although unpublished work shows that the argument can be extended to the group $\Sp(4,\CC)$ of type $\rmB_2 \,/\, \rmC_2$).  
Since we would like to point out the technical point which remains to be resolved for the construction in general, we will begin this discussion in the generality of an arbitrary connected, simply connected, complex semisimple Lie group $G$. 

Let $w \to w'$ be a an edge in the Bruhat graph with $w' = s_\alpha w$ for some \emph{simple} root $\alpha\in\bSigma$.  In this case, we will call it a \emph {simple edge}.    According to Theorem \ref{thm:BGG}, we obtain a $G$-invariant BGG operator
\[
 E_\alpha^n : L^2(\sE_{-w\rho+\rho, -w\rho - \rho}) \to L^2(\sE_{-w'\rho+\rho, -w'\rho - \rho})
\]
for some $n$, and this is a longitudinally elliptic differential operator along the leaves of $\sF_\alpha$.  
We will replace this operator by its bounded operator phase
\begin{equation}
 \label{eq:normalized_BGG_ops}
  T_{w\to w'}:=\ph(E_\alpha^n) : L^2(\sE_{-w\rho+\rho, 0}) \to L^2(\sE_{-w'\rho+\rho, 0}),
\end{equation}
acting between unitary principal series.  We call this the  \emph{normalized BGG operator}.  It belongs to $\Psi^0_\alpha$ and satisfies
\begin{align*}
 T_{w \to w'}^* T_{w\to w'} -I &~\in \Psi^{-1}_\alpha, \\
 [a,T_{w\to w'}] &~\in \Psi^{-1}_\alpha,  &&\forall a\in C^\infty(\sX), \\
 g T_{w\to w'} g^{-1} - T_{w\to w'} &~\in \Psi^{-1}_\alpha,  &&\forall g\in G,
\end{align*}
where $a\in C^\infty(\sX)$ is acting on $\Gamma^{\infty}(\sE_{-w\rho+\rho,0})$ and $\Gamma^{\infty}(\sE_{-w'\rho+\rho,0})$ by multiplication, and $g\in G$ is acting by the unitary principal series representations given in \eqref{eq:normalized_BGG_ops}. 

For the rank-two complex semisimple groups, the system of normalized BGG operators are indicated in Figure \ref{fig:BGG}.
Here $\alpha$ and $\beta$ are the two simple roots, with $\beta$ the longer root in the cases $\rmB_2 \,/\, \rmC_2$ and $\rmG_2$.  
The arrows pointing north-east are $\ph(E_\alpha^n)$ for some $n$, and those pointing south-east are $\ph(E_\beta^n)$ for some $n$.
We have not yet defined the horizontal arrows (associated to the non-simple edges).

\begin{figure}[h]
\scriptsize
\begin{align*}
  &\underline{\rmA_1 \sqcup \rmA_1:} \\
 &\xymatrix@C=1ex@R=1ex{
   & L^2(\sE_{-\alpha,0}) \ar[dr] \ar@{.}[dd]|\bigoplus\\
   L^2(\sE_{0,0}) \ar[ur] \ar[dr] 
   && L^2(\sE_{-\alpha-\beta,0}) \\
   & L^2(\sE_{-\beta,0}) \ar[ur]
 }
\\
 &\underline{\rmA_2:} \\
 &\xymatrix@C=1ex@R=1ex{
   & L^2(\sE_{-\alpha,0}) \ar[ddrr] \ar[rr]\ar@{.}[dd]|\bigoplus && L^2(\sE_{-2\alpha-\beta,0}) \ar[dr]\ar@{.}[dd]|\bigoplus\\
   L^2(\sE_{0,0}) \ar[ur] \ar[dr] 
   &&&& L^2(\sE_{-2\alpha-2\beta,0}) \\
   & L^2(\sE_{-\beta,0}) \ar[uurr] \ar[rr] && L^2(\sE_{-\alpha-2\beta,0}) \ar[ur]
 }
\\
 &\underline{\rmB_2/\rmC_2:} \\
 &\xymatrix@C=1ex@R=1ex{
   & L^2(\sE_{-\alpha,0}) \ar[ddrr] \ar[rr] \ar@{.}[dd]|\bigoplus && L^2(\sE_{-3\alpha-\beta,0}) \ar[ddrr] \ar[rr] \ar@{.}[dd]|\bigoplus && L^2(\sE_{-4\alpha-2\beta,0}) \ar[dr] \ar@{.}[dd]|\bigoplus \\
   L^2(\sE_{0,0}) \ar[ur] \ar[dr] 
   &&&&&& L^2(\sE_{-4\alpha-3\beta,0}) \\
   & L^2(\sE_{-\beta,0}) \ar[uurr] \ar[rr] && L^2(\sE_{-\alpha-2\beta,0}) \ar[uurr] \ar[rr] && L^2(\sE_{-3\alpha-3\beta,0}) \ar[ur]
 }
\\
 &\underline{\rmG_2:} \\
 &\xymatrix@C=0.8ex@R=1ex{
   & L^2(\sE_{-\alpha,0}) \ar[ddrr] \ar[rr] \ar@{.}[dd]|\bigoplus
    && L^2(\sE_{-4\alpha-\beta,0}) \ar[ddrr] \ar[rr] \ar@{.}[dd]|\bigoplus
    && L^2(\sE_{-6\alpha-2\beta,0}) \ar[ddrr] \ar[rr] \ar@{.}[dd]|\bigoplus
    && L^2(\sE_{-9\alpha-4\beta,0}) \ar[ddrr] \ar[rr] \ar@{.}[dd]|\bigoplus
    && L^2(\sE_{-10\alpha-5\beta,0}) \ar[dr] \ar@{.}[dd]|\bigoplus\\
   L^2(\sE_{0,0}) \ar[ur] \ar[dr] 
    &&&&&&&&&& \!\!\!\!L^2(\sE_{-10\alpha-6\beta,0}) \\
   & L^2(\sE_{-\beta,0}) \ar[uurr] \ar[rr] 
    && L^2(\sE_{-\alpha-2\beta,0}) \ar[uurr] \ar[rr] 
    && L^2(\sE_{-4\alpha-4\beta,0}) \ar[uurr] \ar[rr] 
    && L^2(\sE_{-6\alpha-5\beta,0}) \ar[uurr] \ar[rr]
    && L^2(\sE_{-9\alpha-6\beta,0}) \ar[ur]
 }
\end{align*}
\caption{\small The Bernstein-Gelfand-Gelfand complexes of the rank $2$ complex semisimple groups.}
\label{fig:BGG}
\end{figure}

In order to continue we will need to slightly enlarge the $C^*$-algebras $\overline{\Psi^{-\infty}_I}$. 
Let  $\sE = \bigoplus_{w\in W} \sE_{-w\rho+\rho, -w\rho-\rho}$ denote sum of the line bundles in the BGG complex, with $\ZZ/2\ZZ$-grading according to the parity of $\ell(w)$.  Put also $H = L^2(\sE)$.

\begin{definition}
 \label{def:K-A-J}
	For each $ I \subseteq \bSigma$, we define $\KK_I$ to be the hereditary $C^*$-subalgebra of $\LL(H)$ generated by $\overline{\Psi^{-\infty}_I(\sE)}$.  
	We also define $\AA$ to be the simultaneous multiplier algebra of the $\KK_I$,
	\[
	\AA = \{ A \in L^2(H) \st A\KK_I \subseteq \KK_I, \quad \forall I \subseteq \bSigma\}.
	\]
       Finally, we put $\JJ_I = \KK_I \cap \AA$.
\end{definition}

One advantage of the algebras $\KK_I$ over $\overline{\Psi^{-\infty}_I}$ is that they are nested: $\KK_I \subseteq \KK_J$ whenever $J \subseteq I$.  Therefore, the algebras $\JJ_I$ ($I\subseteq \bSigma$) form a nested family of ideals of $\AA$.  

The construction of the BGG class in $K$-homology will take place entirely within the $C^*$-algebra $\AA$.
It is easy to check that multiplication operators by functions on $\sX$ belong to $\AA$, as do the principal series actions of $g\in G$.   We will also need that the normalized BGG operators belong to $\AA$.  The results of \cite{Yuncken:BGG} give this result only for the simply-laced case.

\begin{lemma}
 \label{lem:BBB_in_A}
 Let $G$ be a simply-laced complex semisimple Lie group.  For any $\alpha, \beta\in \bSigma$ we have $\ph(E_\beta).\KK_\alpha \subset \KK_\alpha$.  It follows that the normalized BGG operators belong to $\AA$.
\end{lemma}

\begin{remark}
 \label{rmk:B_2}
 The proof of Lemma \ref{lem:BBB_in_A} depends only on the relative position of two roots $\alpha$ and $\beta$, and so its proof reduces to the rank $2$ groups.  For type $\rmA_1 \sqcup \rmA_1$ the result is elementary.  For $\rmA_2$ it was proven in \cite{Yuncken:BGG} using Gelfand-Tsetlin theory.  In unpublished work we have also generalized this to type $\rmB_2$, which shows that the lemma can be extended to all groups other than $\rmG_2$.  For type $\rmG_2$ it remains unproven, although obviously one expects it to be true as well.
\end{remark}

Next we define normalized BGG operators associated to the non-simple edges.

\begin{lemma}
	\label{lem:shift_trick}
 Suppose $G$ is a group for which Lemma \ref{lem:BBB_in_A} holds.  Then one can associate to each non-simple edge $w \to w'$ in the Bruhat graph a bounded operator 
 \begin{equation}
 \label{eq:T}
  T_{w \to w'} :L^2(\sE_{-w\rho, 0}) \to L^2(\sE_{-w'\rho, 0})
 \end{equation}
 such that the operators $T_{w\to w'}$ form a complex modulo the sum of the ideals $\sum_\alpha \JJ_\alpha$.
\end{lemma}

\begin{proof}
We compare the normalized BGG operators associate to the simple edges, Equation \eqref{eq:normalized_BGG_ops}, with the same operators shifted by $-\rho$:
\begin{equation}
 \label{eq:intertwiners}
  I_{w\to w'}:=\ph(E_\alpha^n) : L^2(\sE_{-w\rho, 0}) \to L^2(\sE_{-w'\rho, 0}).
\end{equation}
The operators $I_{w\to w'}$ are actually intertwiners of irreducible unitary principal series representations (see, \emph{e.g.}\ \cite{Duflo:representations}) so by Schur's Lemma the diagram of the operators $I_{w\to w'}$ will commute on the nose.  

 Using a partition of unity, we can find sections $f_1,\ldots,f_n \in \Gamma^\infty(\sE_{-\rho})$ such that $\sum_{i=1}^n \overline{f_i}\, f_i =1$.  By considering principal symbols, we see that
\begin{equation}
 \label{eq:shift_trick}
 \left(\sum_i \overline{f_i} I_{w\to w'} f_i \right) \;-\; T_{w\to w'} ~\in~  \Psi^{-1}_\alpha \subset \KK_\alpha
\end{equation}
for every edge $w \to  w'$ associated to a simple root $\alpha$.  Let us define the normalized BGG operators associated to non-simple edges $w\to w'$ by 
\[
 T_{w\to w'} := \sum_i \overline{f_i} I_{w\to w'} f_i.
\]
Note that $I_{w\to w'}$ can be written as a composition of intertwiners associated to simple roots (or their inverses).  Using Lemma \ref{lem:BBB_in_A}, we infer that the diagram of normalized BGG operators $T_{w\to w'}$ commutes modulo $\sum_\alpha \JJ_\alpha$.  By introducing signs, as in \cite{BerGelGel}, we obtain a complex modulo $\sum_\alpha \JJ_\alpha$.
\end{proof}

Until this point, the entire construction works for an arbitrary simply laced $G$, and using Remark \ref{rmk:B_2} it can be generalized to any complex semisimple group with no factor of type $\rmG_2$.  The final step, however, is an analogue of the Kasparov product, and this has only been achieved in rank $2$.  We must therefore restrict our attention now to the group $\SL(3,\CC)$.  We also consider the group $\SL(2,\CC) \times \SL(2,\CC)$ in order to make a comparison with the classical Kasparov product.

Theorem \ref{thm:foliations} implies that $\JJ_\alpha. \JJ_\beta \subseteq \KK(H)$ for the rank-two groups.  The Kasparov Technical Theorem yields a pair of bounded self-adjoint operators $N_\alpha, N_\beta \in \LL(H)$ that are diagonal with respect to the direct sum $H = \bigoplus_{w\in W} L^2(\sE_{w\cdot(0,-2\rho)})$,  that commute modulo compacts with the normalized BGG operators, the multiplication action of $C(\sX)$ and the principal series representations of $G$ on each $L^2(\sE_{w\cdot(0,-2\rho)})$, and that satisfy
\begin{align}
 \label{eq:po1}
  N_\alpha \JJ_\alpha \subseteq \KK(H), 
 && N_\beta \JJ_\beta \subseteq \KK(H) \qquad \text{and } 
 && N_\alpha^2 + N_\beta^2 = 1.
\end{align}
We multiply each of the normalized BGG operators by a power of $N_\alpha$ and $N_\beta$ as indicated in Figure \ref{fig:KTT}
\begin{figure}
\footnotesize
\begin{align*}
  &\underline{\rmA_1 \sqcup \rmA_1:} &
 &\hspace{4ex} \xymatrix@C=1.7ex@R=1ex{
   & L^2(\sE_{-\alpha,0}) \ar[dr]^{N_\beta T} \ar@{.}[dd]|\bigoplus\\
   L^2(\sE_{0,0}) \ar[ur]^{N_\alpha T} \ar[dr]_{N_\beta T} 
   && L^2(\sE_{-\alpha-\beta,0}) \\
   & L^2(\sE_{-\beta,0}) \ar[ur]_{N_\alpha T}
 }
\\
 &\underline{\rmA_2:} &
 &\hspace{4ex} \xymatrix@C=1.7ex@R=1ex{
   & L^2(\sE_{-\alpha,0}) \ar[ddrr]|(0.65){N_\beta^2 T} \ar[rr]^{N_\alpha N_\beta T}\ar@{.}[dd]|\bigoplus && L^2(\sE_{-2\alpha-\beta,0}) \ar[dr]^{N_\beta T}\ar@{.}[dd]|\bigoplus\\
   L^2(\sE_{0,0}) \ar[ur]^{N_\alpha T} \ar[dr]_{N_\beta T} 
   &&&& L^2(\sE_{-2\alpha-2\beta,0}) \\
   & L^2(\sE_{-\beta,0}) \ar[uurr]|(0.65){N_\alpha^2 T} \ar[rr]_{N_\alpha N_\beta T} && L^2(\sE_{-\alpha-2\beta,0}) \ar[ur]_{N_\alpha T}
 }
\end{align*}
\caption{\small Modification of the normalized BGG operators by the Kasparov Technical Theorem.}
\label{fig:KTT}
\end{figure}
where the symbol $T$ is shorthand for the normalized BGG operator $T_{w\to w'}$ associated to the given edge.  Finally, we define $F$ to be the sum of all these operators and their adjoints.  A direct check using the properties of $N_\alpha$, $N_\beta$ and the normalized BGG operators shows that $F^2 \equiv 1$ modulo $\KK(H)$.

We arrive at the following result.

\begin{theorem}
 \label{thm:KK_BGG}
 Let $G=\SL(2,\CC) \times \SL(2,\CC)$ or $\SL(3,\CC)$ and let $\sX$ be the full flag manifold of $G$.  The operator $F$ described above defines a $G$-equivariant $K$-homology class for $\sX$, which we denote by $\BGG \in KK^G(C(\sX),\CC)$.  Its image under the forgetful map to $KK^K(\CC,\CC)$ is the class of the trivial representation, and its image in $KK^G(\CC,\CC)$ is Kasparov's $\gamma$ element.
\end{theorem}

For the group $G=\SL(2,\CC) \times \SL(2,\CC)$, the flag manifold is $\sX = \CP^1 \times \CP^1$ and the four induced bundles in Figure \ref{fig:KTT} are isomorphic to
 \[
  L^2(\sE_{-i\alpha-j\beta}) \cong L^2 (\Lambda^{0,i}\CP^1 \,\otimes \; \Lambda^{0,j}\CP^1), \qquad i,j\in\{0,1\}. 
 \]
In other words, Figure \ref{fig:KTT} for $\rmA_1 \sqcup \rmA_1$ depicts precisely the exterior Kasparov product of two copies of the Dolbeault-Dirac class for $\CP^1$.  Thus the BGG class for $\SL(3,\CC)$ is literally a generalization of the Kasparov product.

If we allow ourselves to use Lemma \ref{lem:BBB_in_A} for the group of type $\rmB_2$, as indicated in Remark \ref{rmk:B_2}, then we can obtain Theorem \ref{thm:KK_BGG} also for the group $\Sp(4,\CC)$ of type $\rmB_2 \, / \, \rmC_2$ by modifying the normalized BGG complex as in Figure \ref{fig:KTT2}.  

\begin{figure}[h]
\footnotesize
\begin{align*}
 &\underline{\rmB_2/\rmC_2:} &
 & \xymatrix@C=1.7ex@R=1ex{
   & L^2(\sE_{-\alpha,0}) \ar[ddrr]|(0.65){N_\beta^2 T} \ar[rr]^{N_\alpha N_\beta T} \ar@{.}[dd]|\bigoplus && L^2(\sE_{-3\alpha-\beta,0}) \ar[ddrr]|(0.65){N_\beta^2 T} \ar[rr]^{N_\alpha N_\beta T} \ar@{.}[dd]|\bigoplus && L^2(\sE_{-4\alpha-2\beta,0}) \ar[dr]^{N_\beta T} \ar@{.}[dd]|\bigoplus \\
   L^2(\sE_{0,0}) \ar[ur]^{N_\alpha T} \ar[dr]_{N_\beta T} 
   &&&&&& L^2(\sE_{-4\alpha-3\beta,0}) \\
   & L^2(\sE_{-\beta,0}) \ar[uurr]|(0.65){N_\alpha^2 T} \ar[rr]_{N_\alpha N_\beta T} && L^2(\sE_{-\alpha-2\beta,0}) \ar[uurr]|(0.65){N_\alpha^2 T} \ar[rr]_{N_\alpha N_\beta T} && L^2(\sE_{-3\alpha-3\beta,0}) \ar[ur]_{N_\alpha T}
}
\end{align*}
\caption{\small Modification of the normalized BGG operators in type $\rmB_2 \,/\, \rmC_2$.}
\label{fig:KTT2}
\end{figure}

\chapter{Pseudodifferential operators on quantum flag manifolds}
\label{sec:CQG}

Quantum groups and Connes' noncommutative geometry arose independently, but it was soon realized that algebras of functions on quantized Lie groups and their homogeneous spaces should, morally, be examples of noncommutative manifolds.  Unfortunately, making this precise has turned out to be much more difficult than expected.

There is no definitive list of axioms of a noncommutative manifold, rather a list of desirable properties which should be adapted to suit the various examples.  The basic object, however, is always a \emph{spectral triple} $(\sA,H,D)$ where $\sA$ is a $*$-algebra of bounded operators on the Hilbert space $H$, and $D$ is an unbounded self-adjoint operator on $H$ satisfying axioms designed to mimic the properties of a Dirac-type operator.  This needs to be supplemented, at minimum, with some version of the \emph{regularity} axiom \cite{Connes:NCG, ConMos:local_index_formula}.  See below for definitions.  

In practice, regularity amounts to the existence of a pseudodifferential calculus, in which $D$ is elliptic.  This was made into an explicit theorem by Higson \cite{Higson:local_index_formula} and Uuye \cite{Uuye:PsiDOs}.  Not only does regularity imply the existence of an abstract pseudodifferential calculus, but proving the regularity axiom is typically achieved by exhibiting that pseudodifferential calculus.

In this chapter, we shall discuss the noncommutative geometry of two classes of quantum homogeneous spaces by studying analogues of pseudodifferential operators upon them.

We first look at the quantum projective spaces $\CP^n_q$, chosen as a case where things work relatively easily.  Quantum projective spaces admit spectral triples which are $q$-deformations of the classical Dirac-type operators (\cite{Krahmer:Dirac},\cite{DAnDab}).  They fail to be regular, but only very mildly.  In the case of the Podle\'{s} sphere $\CP^1_q$, Neshveyev and Tuset produced a twisted algebra of pseudodifferential operators and thus proved a local index formula for the Dabrowski-Sitarz spectral triple \cite{NesTus:local_index_formula}.   We will present some recent work with M.~Matassa giving the first steps towards a similar result for general $\CP^n_q$.  We will generalize the pseudodifferential calculus of Connes and Moscovici to twisted spectral triples, in a sense slightly more general than \cite{ConMos:twisted, Moscovici:twisted}.  With this definition, the $q$-deformed Dirac operators of \cite{KraTuc} yield twisted regular spectral triples.

In a second example, we consider the quantized full flag manifold $\sX_q$ of $K_q = \SU_q(3)$.  Here, the usual axioms of a noncommutative manifold fail more spectacularly.  We don't even have candidates for Dirac operators on $\sX_q$ which are analogous the Dabrowski-Sitarz operator on $\CP^1$, and there are good reasons to be pessimistic about the existence of such operators (although see \cite{NesTus:isospectral} for a very different approach to defining spectral triples on quantum groups).

On the other hand, there is a $q$-deformation of the Bernstein-Gelfand-Gelfand complex on $\sX_q$.  We shall describe the results of \cite{VoiYun:SUq3}, which show that the rudimentary multifiltered pseudodifferential theory of Chapter \ref{sec:flags} survives the quantization process.  Once again, this is far from being as powerful as a full pseudodifferential calculus, but it still suffices to obtain some useful results.  As in the classical case, we can obtain an $G_q$-equivariant Bernstein-Gelfand-Gelfand class in $KK^{G_q}(C(\sX_q),\CC)$, where $G_q$ denotes the Drinfeld double $K_q \!\bowtie\! \hat{K}_q$.  As a result, following earlier work of Meyer, Nest and Voigt \cite{MeyNes:localisation, MeyNes:coactions, NesVoi, Voigt:free_orthogonal},
we obtain equivariant Poincar\'e duality for $\sX_q$ and a proof of the Baum-Connes Conjecture for the discrete quantum dual of $\SU_q(3)$.

\section{Notation and conventions}
\label{sec:quantum_notation}

We continue to use the notation of Section \ref{sec:flag_notation} for complex semisimple Lie groups and their Lie algebras. For their quantizations, we follow the conventions of \cite{VoiYun:CQG}, which we rapidly recap here.   Unless otherwise stated, we assume $0<q<1$.

\begin{itemize}
 \item $\Uqg$ is the Hopf algebra with 
  \begin{itemize}
    \item generators: $E_\alpha$, $F_\alpha$ (for $\alpha\in\bSigma$) and $K_\lambda$ ($\lambda \in \bP$),
    \item algebra relations:
 \begin{align*}
  &K_0 = 1, && K_\lambda K_\mu = K_{\lambda+\mu}, \\
  &K_\lambda E_\alpha K_\lambda^{-1} = q^{(\lambda,\alpha)} E_\alpha, 
  &&K_\lambda F_\alpha K_\lambda^{-1} = q^{-(\lambda,\alpha)} E_\alpha, \\
  &[E_\alpha, F_\beta] = \delta_{\alpha,\beta} \frac{K_\alpha - K_\alpha^{-1}}{q-q^{-1}},
 \end{align*}
 and the quantum Serre relations, which we shall not write out here,
    \item coalgebra relations:
 \begin{align*}
  &\hat\Delta E_\alpha = E_\alpha \otimes K_\alpha + 1 \otimes E_\alpha, &&
  \hat\Delta F_\alpha = F_\alpha \otimes 1 + K_\alpha^{-1} \otimes F_\alpha,  \\
  &\hat\Delta K_\lambda = K_\lambda \otimes K_\lambda, \quad \\
  &\hat\epsilon(E_\alpha) = \hat\epsilon(F_\alpha) = 0, && \hat\epsilon(K_\lambda) = 1,
 \end{align*}
    \item antipode:\qquad $
  \hat{S}(K_\lambda) = K_{-\lambda},  \quad 
  \hat{S}(E_\alpha) = - E_\alpha K_{\alpha}^{-1},  \quad
  \hat{S}(F_\alpha) = -K_\alpha F_\alpha.
 $
 \end{itemize}

 \item $\UqRk$ denotes $\Uqg$ with the $*$-structure 
  \[
   E_\alpha^* = K_\alpha F_\alpha, \qquad  F_\alpha^* = E_\alpha K_\alpha^{-1}, \qquad K_\alpha^*=K_\alpha.
  \]

 \item $\sA(K_q)$ denotes the algebra of matrix coefficients of finite dimensional $\UqRk$-modules of type $1$ (meaning that the $K_\lambda$ act as positive operators), with Hopf $*$-algebra structure dual to $\UqRk^\cop$, \emph{i.e.}, for all $X,Y \in \UqRk$, $a,b\in\sA(K_q)$,
 \begin{align*}
  (X,ab) &= (X_{(1)},b)(X_{(2)},a),   & (XY,a) = (X,a_{(1)})(Y,a_{(2)}), \\
    (X,a^*) &= \overline{(\hat{S}^{-1}(X)^*, a)}.
 \end{align*}

 \item The correspondence between left $\sA(K_q)$-comodules and left $\UqRk$-modules (of type $1$) is made by equipping a left $\sA(K_q)$-comodule $V$ with the left $\UqRk$-action
 \[
  X. v := (S(X), v_{(-1)}) v_{(0)} \qquad  \forall X\in \UqRk, ~v\in V.
 \]
 
 \item For $X \in \UqRk$, $a\in \sA(K_q)$, we will write
 \begin{equation}
  \label{eq:hit}
  X\hit a = a_{(1)} (X,a_{(2)}),
 \end{equation}
 for the right regular representation of $\sU(\lie{k})$ on $C^\infty(K)$.

 \item $\phi$ denotes the bi-invariant Haar state on $\sA(K_q)$, $L^2(K_q)$ is the corresponding GNS space, and $C(K_q)$ is the $C^*$-closure of the GNS representation of $\sA(K_q)$.

\end{itemize}

If $L_q$ is a quantum subgroup of $K_q$, defined via a surjective morphism $\res_{L_q}:\sA(K_q) \onto \sA(L_q)$, then the quantum homogeneous space $K_q/L_q$ is defined by its {algebra of polynomial functions},
\[
 \sA(K_q/L_q) = \{ a\in \sA(K_q) \st (\id\otimes\res_{L_q})\Delta a = a \otimes 1 \}.
\]
The closures of $\sA(K_q/L_q)$ in $L^2(K_q)$ and $C(K_q)$ are denoted by $L^2(K_q/L_q)$ and $C(K_q/L_q)$, respectively.

If $V$ is a finite-dimensional unitary left $\sA(L_q)$-comodule, with coaction $\alpha: V \to \sA(L_q)\otimes V$, then the \emph{induced bundle} $K_q\times_{L_q} V$ is defined via its space of polynomial sections:
\begin{align*}
 \sA(K_q \times_{L_q} V) &= \{ \xi \in \sA(K_q) \otimes V \st 
   (\id \otimes \res_{L_q} \otimes \id)(\Delta \otimes \id) \xi = (\id \otimes \alpha)\xi \}.
\end{align*}
Equivalently, an element $ \xi = \sum_i a_i \otimes v_i \in \sA(K_q) \otimes V$ belongs to $\sA(K_q \times_{L_q} V)$ if and only if
\begin{equation}
\label{eq:invariant}
 \sum_i (X_{(1)}\hit a_i) \otimes (X_{(2)}. v_i) = \hat\epsilon(X) \sum_i a_i \otimes v_i
  \qquad\qquad \forall X \in \UqRl.
\end{equation}
The space $\sA(K_q\times_{L_q} V)$ is a left $\sA(K_q/L_q)$-module by left multiplication on the first leg, and a left $\sA(K_q)$-comodule by comultiplication on the first leg.  It's closure in $L^2(K_q) \otimes V$ is denoted $L^2(K_q \times_{L_q} V)$.

\section{Pseudodifferential calculus and twisted spectral triples}
\label{sec:twisted}

\subsection{Twisted regularity}

In this section and the next, a \emph{twisting} of an algebra $\sA$ will mean a linear automorphism $\theta:\sA \to \sA$ (not necessarily an algebra automorphism).  We will use the notation
\[
  [x,y]_\theta = xy-\theta(y)x
\]
for the \emph{twisted commutator} of $x$ and $y$.

\begin{definition}
	\textup{(\cite{ConMos:twisted,Moscovici:twisted,MatYun:PsiDOs})}
	A \emph{(unital) twisted spectral triple} (called a \emph{type III spectral triple} in \cite{ConMos:twisted}) is $(\sA,H,D,\theta)$ where
	\begin{itemize}
		\item $H$ is a Hilbert space,
		\item $\sA$ is a unital $*$-algebra, represented as bounded operators on $H$, with twisting $\theta$,
		\item $D$ is an unbounded self-adjoint operator on $H$,
	\end{itemize}
  such that
  \begin{enumerate}
  	\item $D$ has compact resolvent,
  	\item $[D,a]_\theta$ is densely defined and bounded for all $a\in\sA$.
  \end{enumerate}
  It is \emph{regular} if $\theta$ can be extended to a $*$-algebra $\sB\subseteq\LL(H)$ of bounded operators containing both $\sA$ and $[D,\sA]_\theta$ such that the twisted derivation
  \[
   \delta_\theta = [|D|,\slot]
  \]
  preserves $\sB$.
\end{definition}

Our definition of a twisted spectral triple is weaker than the definitions of Connes and Moscovici, thanks our allowing $\theta$ to be merely a linear isomorphism, not an algebra automorphism.  This relaxation of the twisting will be necessary for application to quantum homogeneous spaces.

Our main goal will be to prove a twisted analogue of the results of \cite{Higson:local_index_formula, Uuye:PsiDOs} which show the equivalence of regularity and the existence of a pseudodifferential calculus.

\subsection{Differential and pseudodifferential operators}

The following definitions are standard. 
\begin{itemize}
	\item $\Delta = D^2+1$ is referred to as the \emph{Laplace operator},
	\item $H^\infty = \bigcap_{n\in\NN} \Dom(\Delta^n)$, is the space of \emph{smooth elements} in $H$,
	\item The \emph{$s$-Sobolev space} $H^s$ (for $s\in\RR$) is the completion of $H^\infty$ with respect to the norm
	\[
	 \| v \|_s = \|\Delta^{\frac{s}{2}}v\|,
	\]
	\item  $\Op^t$ (for $t\in \RR$) is the set of  operators $T$ on $H^\infty$ which extend to bounded operators $T:H^{s+t} \to H^s$ for every $s\in\RR$; these are called the \emph{operators of analytic order $t$}.
\end{itemize}

\begin{remark}
	The Laplace operator $\Delta$ is intended as an invertible replacement of $D^2$.  Note that $|D| - \Delta^{\frac12} \leq \Delta^{-\frac12} \in Op^{-\frac12}$.  This will mean that $|D|$ and $\Delta^\half$ are essentially interchangeable in the analysis that follows, and we will occasionally do so without further comment.
\end{remark}

\begin{definition}
	\label{def:DOs}
	An \emph{algebra of differential operators adapted to a twisted spectral triple $(\sA,H,D,\theta)$} is an $\NN$-filtered algebra $\DO^\bullet$ of operators on $H^\infty$ which contains both $\sA$ and $[D,\sA]_\theta$ in $\DO^0$, is equipped with an extension of $\theta$, and satisfies:
		\begin{enumerate}
			\item \emph{Elliptic estimates}: for every $X\in\DO^n$, there exists $C>0$ such that 
  \[
    \|Xv\|_H \leq C\|\Delta^{\frac{n}{2}}\|_H \qquad \text{ for all } v\in H^\infty, 
  \]

			\item $[\Delta,\DO^n]_{\theta^2} \subseteq \DO^{n+1}$ for all $n$.
		\end{enumerate}
\end{definition}

\begin{definition}
  An \emph{algebra of pseudodifferential operators adapted to $(\sA,H,D,\theta)$} is an $\RR$-filtered subalgebra $\Psi^\bullet$ of $\Op^\bullet$ containing both $\sA$ and $[D,\sA]_\theta$ in $\Psi^0$, equipped with a $1$-parameter group $(\Theta^z)_{z\in\CC}$ of twistings, and satisfying
  \begin{enumerate}
  	\item $\theta(a)-\Theta^1(a) \in \Psi^{-1}$ for all $a\in\sA$,
  	\item $\Delta^z \Psi^\bullet \subseteq \Psi^\bullet$ and $\Psi^\bullet\Delta^z \subseteq \Psi^\bullet$ for all $z\in\CC$,
  	\item $[\Delta^{\frac{z}{2}},\Psi^s]_{\Theta^z} \subseteq \Psi^{\Re(z)+s-1}$ for all $z\in\CC$, $s\in\RR$.
  \end{enumerate}
\end{definition}

\begin{theorem}\textup{(\cite{MatYun:PsiDOs})}
  \label{thm:regularity}
	Let $(\sA,H,D,\theta)$ be a twisted spectral triple.  
	\begin{enumerate}
		\item If $(\sA, H, D,\theta)$ is regular, then there exist algebras of differential operators and pseudodifferential operators adapted to it.
		\item Conversely, if there exists an algebra of differential operators $\DO$ adapted to $(\sA,H,D,\theta)$, with $\theta$ \emph{diagonalizable}---meaning that $\DO$ is the direct sum of its $\theta$-eigenspaces---then $(\sA,H,D,\theta)$ is regular.
	\end{enumerate}
\end{theorem}

\begin{proof}
 The first statement is straightforward.  Using $\delta_\theta$-stability, one can prove that the algebra $\sB$ from the definition of regularity lies in $\Op^0$.  Then we can define $\Psi^\bullet$ as the subalgebra of $\Op^\bullet$ generated by $\sB$ and $\Delta^z$ for all $z\in\CC$, equipped with the family of twistings
 \[
   \Theta^z(X) = \Delta^{z} X \Delta^{-z}, \qquad X\in \Psi, z\in\CC.
 \]
 This yields an algebra of pseudodifferential operators adapted to $(\sA,H,D,\theta)$.  For the algebra of differential operators, put $\DO^m = \Psi^m$.  	We have $\Theta^1(b) - \theta(b) \in \Psi^{-1}$ for all $b\in \sB$, and consequently we can extend $\theta$ to a perturbation of $\Theta^1$ on $\DO$ which satisfies the axioms of an algebra of differential operators. 

  The second statement is more profound.  Suppose $\DO^\bullet$ is an algebra of differential operators adapted to $(\sA,H,D,\theta)$, with $\theta$ diagonalizable.  It suffices to define an adapted algebra of pseudodifferential operators, since then we can use $\sB=\Psi^0$ to obtain regularity.   To define the pseudodifferential operators, we start with the following definition.

\begin{definition}
  Let $T$ be an operator on $H^\infty$.  We say that \emph{$T$ has an asymptotic expansion $T \sim \sum_{i=1}^\infty T_i$}, where each $T_i\in\Op^\bullet$, if for any $r\in \RR$ we have $T-\sum_{i=1}^N T_i \in \Op^{-r}$ for all sufficiently large $N$.
\end{definition}

  One defines a \emph{basic pseudodifferential operator} to be an operator $T$ on $H^\infty$ that admits an asymptotic expansion of the form
  \[
   T \sim \sum_{k=0}^N X_k \Delta^{\frac{z}{2}-k},
  \]
  with $X_k \in \DO$.
  We then define $\Psi$ as the space of finite linear combinations of basic pseudodifferential operators.

  It is nontrivial to prove that $\Psi$ is an algebra.  The key point is to prove that if $X\in\DO$ we have $\Delta^z X \in \Psi$. The standard trick here is use the Cauchy Integral Formula to write, when $\Re(z) <0$,
\[
 \Delta^z  = \oint_\Gamma \lambda^z(\lambda-\Delta)^{-1} \,d\lambda,
\]
where $\Gamma$ is a vertical line separating the spectrum of $\Delta$ from $0$.  In the untwisted case, one can then repeatedly apply the commutator formula
\[
 [(\lambda - \Delta)^{-1} , Y] = (\lambda-\Delta)^{-1} [\Delta,Y] (\lambda-\Delta)^{-1}, \qquad (Y\in\DO)
\]
to develop the product $\Delta^z X$ as an asymptotic expansion.

In the twisted case, the commutators with the resolvent are more delicate, and it's here we take advantage of the diagonalizability of $\theta$.  Let us write $W= \Sp(\theta)$ for the set of eigenvalues of $\theta$. 

We write $\nabla(X) = [\Delta,X]_{\theta^2}$ for the twisted commutator of $X \in \DO$.  
We will define a decomposition of the iterated twisted commutator $ \nabla^n(X)$ indexed by sequences $\bfmu = (\mu_0,\ldots,\mu_n) \in W^n$, as follows.  If $n=0$ then $\nabla^{(\mu_0)}(X)$ denotes the $\mu_0$-eigencomponent of $X$. If $n>0$ we define recursively 
\[
 \nabla^\bfmu(X) = \nabla(\nabla^{(\mu_0,\ldots,\mu_{n-1})}(X))^{\mu_n}.
\]
Thus 
\[
 \nabla^n(X) = \!\! \sum_{\bfmu\in W^{n+1}} \!\!\nabla^\bfmu(X).
\]

\begin{lemma}
 For any $X\in\DO$ and $z\in\CC$, we have the asymptotic expansion
 \[
  \Delta^z X \sim \sum_{n=0}^\infty \sum_{\bfmu\in W^{n+1}} {z \choose n}_{\!\!\bfmu} \nabla^\bfmu(X) \Delta^{z-k},
 \]
 where the ${z \choose n}_{\bfmu}$ are the quantized binomial coefficients defined in \cite[Appendix A]{MatYun:PsiDOs}.  The $n=0$ term is equal to $\theta^z (X) \Delta^{z-k}$.
\end{lemma}

From this, we obtain that $\Psi^\bullet$ is an algebra.  We can define a one parameter family of twistings on $\Psi$ by $\Theta^z(T) = \Delta^{z/2}T\Delta^{-z/2}$, and it is easy to check that $\Psi^\bullet$ satisfies all the axioms of an algebra of pseudodifferential operators.
\end{proof}

\subsection{Pseudodifferential calculus on quantum projective spaces}

The framework above was designed with the following application in mind.   Let $G=\SL(n+1,\CC)$ with simple roots $\{\alpha_1,\ldots,\alpha_n\}$.  Let $\bar{P} = \bar{P}_{\{\alpha_1,\ldots,\alpha_{n-1}\}}$ be the lower parabolic subgroup corresponding to $I=\{\alpha_1,\ldots,\alpha_{n-1}\}\subseteq \bSigma$---see Section \ref{sec:flag_notation}---so that $G/\bar{P} = \CP^n$.  Let $K=SU(n+1)$ and $L=\bar{P}\cap K = S(U(n)\times U(1))$, and note that $K/L = \CP^n$.  

We denote the nilpotent radical of $\bar{\lie{p}}$ by $\bar{\lie{n}}$, and observe that $\bar{\lie{n}} \cong (\lie{g}/ \bar{\lie{p}})^*$, which is an irreducible $L$-module via the coadjoint action.  The anti-holomorphic exterior bundle of $\CP^1$ is isomorphic to $K\times_L \Lambda(\bar{\lie{n}})$

Now let $\sA(K_q)$, $\sA(L_q)$, $\sA(\CP_q^n)$ be the algebras of polynomial functions on $K_q$, $L_q$ and $\CP_q^n$.  Equip $\bar{\lie{n}}$ with the left $\sU_q(\lie{l})$ representation which corresponds to the $L$-module structure above.  Kr\"ahmer and Tucker-Simmons \cite{KraTuc} showed that one can define a $q$-deformation $\Lambda_q(\bar{\lie{n}})$ of the exterior algebra of $\bar{\lie{n}}$.  They use this to define a natural $q$-deformation of the Dirac-Dolbeault operator $D = \bar\partial + \bar\delta^*$ on the sections of the quantum exterior bundle $\Omega_q := K_q \otimes_{L_q} \Lambda_q(\bar{\lie{n}})$.

\begin{theorem}
	\label{thm:D_regular}
 Let $\sA = \sA(\CP^n_q)$, $H=L^2(\Omega_q)$ and $D$ be the Dirac-Dolbeault operator of \cite{KraTuc}. The twisted spectral triple $(\sA, H, D, \id)$ is regular.
\end{theorem}

\begin{remark}
 The fact that $(\sA,H,D)$ has compact resolvent was proved in \cite{DAnDab}.
\end{remark}

Note that the twisting is trivial on the algebra $\sA(\CP^n_q)$, but the spectral triple is not regular in the untwisted sense.  The associated algebras of differential and pseudodifferential operators will be twisted, as seen for instance in \cite{NesTus:local_index_formula}.
In the rest of this section, we will sketch the proof of Theorem \ref{thm:D_regular}.

 We have to define a twisted algebra of differential operators.  We begin by defining the algebra of differential operators on $K_q$ (with polynomial coefficients) as the smash product
 \[
  \DO(K_q)  = \sA(K_q) \# \UqRk.
 \]
 That is, $\DO(K_q) = \sA(K_q) \otimes \UqRk$ as a space, with algebra structure given by the usual products on $\sA(K_q)$ and $\UqRk$ and the commutation relation
 \[
  Xa = (X_{(2)}\hit a) X_{(1)}, \qquad\qquad \forall X\in\UqRk,~a\in\sA(K_q).
 \]
 These relations are defined so that the multiplication action of $\sA(K_q)$ and the $\rightharpoonup$ action \eqref{eq:hit} of $\UqRk$ define an action of $\DO(K_q)$ on $\sA(K_q)$.
 
 We define an algebra filtration on $\UqRk$ by declaring
 \[
  \UqRl \subset \sU_q^0(\lie{k}), \qquad E_{\alpha_n}, F_{\alpha_n} \in \sU_q^1(\lie{k}).
 \]
 This is a Hopf algebra filtration, meaning that $\hat\Delta \sU_q^m(\lie{k}) \subseteq \sum_{i=0}^m \sU_q^i(\lie{k}) \otimes \sU_q^{m-i}(\lie{k})$.  We can therefore extend the filtration to $\DO(K_q)$ by putting
 $
  \DO^m(K_q) = \sA(K_q) \# \sU^m_q(\lie{k}).
 $

 Now let $V$ be a finite dimensional $L_q$-module and let $\sE = K_q\times_{L_q} V$ be the associated induced bundle over $\CP^n_q$.  Define an adjoint action of $\UqRl$ on $\DO(K_q) \otimes \End(V)$ by
 \begin{multline*}
  \ad(X) (A \otimes T) = X_{(1)} A S(X_{(4)}) \otimes X_{(2)} T S(X_{(3)}), \\
 \end{multline*}
 for $X\in\UqRl$, $A \otimes T \in \DO(K_q) \otimes \End(V)$.  From Equation \eqref{eq:invariant} it follows that the invariant subalgebra
 \begin{multline*}
  \DO(\sE) := (\DO(K_q) \otimes \End(V))^{\ad(\UqRl)} \\
    = \{ P \in \DO(K_q) \otimes \End(V) \st \ad(X) P = \hat\epsilon(X) P \quad \forall X \in \UqRl \}
 \end{multline*}
 preserves the section space $\sA(\sE) \subseteq \sA(K_q) \otimes V$.  We refer to $\DO(\sE)$ as the algebra of differential operators on $\sE$, with filtration inherited from $\DO(K_q)$. 

The Dirac-Dolbeault operator $D$ of \cite{KraTuc} belongs to $\DO^1(\Omega_q)$.  D'Andrea-Dabrowski proved that it satisfies the following Parthasarathy formula.

\begin{proposition} \textup{\cite{DAnDab}}
 There is a central element $\sC \in \sU^2_q(\lie{k})$ such that
 \[
   D^2 - \sC \otimes I_{\Lambda_q(\bar{\lie{n}})} \quad  \in \DO^1(\Omega_q).
 \]
\end{proposition}
 
The element $\sC$ is called the \emph{(order $2$) Casimir element}.  D'Andrea-Dabrowski show that, as an unbounded operator on $L^2(\Omega_q)$, $\sC \otimes I$ has compact resolvent, so $(\sA,H,D)$ is a spectral triple.  From these calculations, we also see that $E_{\alpha_n}^*E_{\alpha_n}, F_{\alpha_n}^* F_{\alpha_n} \leq \sC$ as unbounded operators, from which we can deduce the elliptic estimates of Definition \ref{def:DOs}.

The coproduct of $\sC$ has the particularly nice form
\[
  \hat\Delta(\sC) = K_{\omega_n}^2 \otimes \sC + Q + \sC \otimes K_{\omega_n}^{-2} 
\]
for some $Q \in \sU^1_q(\lie{k}) \otimes \sU^1_q(\lie{k})$, where $\omega_n\in \bP$ is the fundamental weight satisfying $(\omega_n,\alpha_i) = \delta_{in}$ for all $i$.  It follows that, for any $a\in\sA(K_q)$,
\[
  \sC a - (K_{\omega_n}^{-2}\hit a) \sC \in \DO^1(K_q)
\]
and for any $X\in\UqRk$,
\[
  \sC X - X \sC = 0.
\]
We are led to define the following twisting.

\begin{definition}
 	Define a twisting $\theta$ on $\DO(K_q)\otimes\End(\Lambda_q(\bar{\lie{n}}))$ via the formula
	\[
	 \theta(aX\otimes T) = (K_{\omega_n}^{-2}\hit a)X \otimes T
	\]
	for all $a\in\sA(K_q)$, $X\in\UqRk$, $T\in\End(\Lambda_q(\bar{\lie{n}}))$.
\end{definition}

Note also that $K_{\omega_n}^2$ commutes with all elements of $\sU_q(\lie{l})$.  
We obtain the following.

\begin{lemma}
 The twisting $\theta$ on $\DO(K_q)\otimes\End(\Lambda_q(\bar{\lie{n}}))$ preserves the subspace $\DO(\Omega_q)$.  With this twisting, $\DO(\Omega_q)$ is a twisted algebra of differential operators adapted to $(\sA,H,D,\id)$.  
\end{lemma}

This completes the proof of Theorem \ref{thm:D_regular}.


\section[Pseudodifferential operators on the full flag manifold of $\SU_q(n)$]{Pseudodifferential operators on the full flag\linebreak manifold of $\SU_q(n)$}
\label{sec:q-flags}

We now pass to our final example, the $q$-deformation of the full flag manifold of $\SU_q(n)$.  Here we will emulate the algebras of longitudinal pseudodifferential operators defined in Section \ref{sec:flags}.

On the full quantum flag manifold of a quantized compact semisimple Lie group $K_q$, there are two reasons why the Bernstein-Gelfand-Gelfand complex seems to be the appropriate point of departure for equivariant index theory.  The first, as already mentioned, is that there is no reasonable candidate for a Dirac operator analogous to those just discussed for $\CP^n_q$.  

More profoundly, as Nest and Voigt pointed out in \cite{NesVoi}, the Kasparov product is not well-defined in $K_q$-equivariant $KK$-theory.  In order to define a product of two $C(K_q)$-comodule algebras, one of them needs to be equipped with a Yetter-Drinfeld module structure  (see also \cite{Vaes:induction}), or equivalently, a representation of the Drinfeld double $K_q \!\bowtie\! \hat{K}_q$.  The Drinfeld double of a $q$-deformed complex semisimple Lie group $K_q$ is often referred to as its \emph{complexification} $G_q$ because it behaves as a $q$-deformation of the complexification $G$ of $K$.  Thus in order to have a $KK$-product, we need to work in $G_q$-equivariant $KK$-theory.  In light of the results of Chapter \ref{sec:flags}, this again suggests the BGG complex.

In the remaining sections, we will show that the BGG class in the $\SL(3,\CC)$-equivariant $K$-homology of the full flag manifold $\sX$ does indeed admit a $q$-deformation to a $\SL_q(3,\CC)$-equivariant $K$-homology class for $\sX_q$.  Again, we will not define a full pseudodifferential calculus, but merely a family of algebras of "operators of negative order" associated to each of the canonical fibrations of the quantum flag manifold.  This time, we will rely upon harmonic analysis of the compact quantum subgroups associated to the various fibrations, rather than the nilpotent subgroups.


\subsection{Algebras of longitudinal pseudodifferential operators on quantum flag manifolds}

We will continue with the notation from Sections \ref{sec:flag_notation} and \ref{sec:quantum_notation}, although we will specialize to the case $K = \SU(n+1)$.  We introduce the following notation for quantum flag manifolds.  Here, $I,J\subseteq\bSigma$ are sets of simple roots.
\begin{itemize}
	\item $K_I = K \cap P_I$ denotes the compact part of the parabolic subgroup $P_I$, and $K_{I,q}$ its $q$ deformation.  The latter is a quantum subgroup of $K_q$ via a surjective morphism $\res_{K_{I,q}}:\sA(K_q) \onto \sA(K_{I,q})$.  In particular, $K_{\bSigma,q} = K_q$ and $K_{\emptyset,q} = T$, the classical torus subgroup.
	\item $\sX_{I,q} = K_q / K_{I,q}$ is the quantized partial flag manifold of $K_q$.  In particular, we write $\sX_q = \sX_{\emptyset,q} = K_q/T$ for the full quantum flag  manifold.
	\item If $\mu\in\bP$ is an integral weight, we write $e^\mu \in C(T)$ for the corresponding unitary character of $T$, and $\CC_\mu$ for the associated one-dimensional left corepresentation of $C(T)$, given by $z \mapsto e^\mu\otimes z$ for all $z\in\CC$.
	\item $\sE_\mu = K_q\times_T \CC_\mu$ is the induced line bundle over $\sX_q$, with section space
        \[
         \sA(\sE_\mu) = \{\xi\in \sA(\sX_q) \st (\id\otimes\res_T)\Delta\xi = \xi\otimes e^\mu\}.
        \]
        The Hilbert completion with respect to the Haar state on $\sA(K_q)$ is denoted $L^2(\sE_\mu)$.
\end{itemize}

To characterize the negative order pseudodifferential operators, we are obliged to use harmonic analysis with respect to the subgroups $K_{I,q}$.  We introduce further notation:
\begin{itemize}
  \item $\Irr(K_{I,q})$ is the set of finite dimensional irreducible $K_{I,q}$-representations, up to equivalence.  

  \item If $\sigma\in\Irr(K_{I,q})$, $V^\sigma$ denotes the vector space underlying $\sigma$.

  \item The algebras of functions and compactly supported functions on the discrete dual of $K_{I,q}$ are, respectively, 
  \[
    \sA(\hat{K}_{I,q}) =  \!\! \prod_{\sigma\in\Irr(K_{I,q})} \!\! \End(V^\sigma),
    \qquad\qquad
    \sAc(\hat{K}_{I,q}) =  \!\! \bigoplus_{\sigma\in\Irr(K_{I,q})} \!\! \End(V^\sigma).
  \]
  Note that there are canonical embeddings $\sA(\hat{K}_{I,q}) \into \sA(\hat{K}_q)$.

  \item For any $I\subset\bSigma$ and any $\sigma \in \Irr(K_{I,q})$, $p_\sigma$ denotes the central projection in $\sAc(\hat{K}_{I,q})$, which acts on any unitary $K_q$-representation as the projection onto the $\sigma$-isotypical subspace.   
  
  \item In the particular case $I=\emptyset$, we identify $e^\mu\in\Irr(T)$ with $\mu\in\bP$.  Then $p_\mu$ acts as the projection onto the $\mu$-weight space of any unitary $K_q$-representation.
\end{itemize}

In particular, we have
\[
 L^2(\sE_\mu) = p_\mu \hit L^2(K_q),
\]
where $\hit$ denotes the obvious extension of the \emph{right regular representation} \eqref{eq:hit} to $\sA(\hat{K}_q)$.

The following lemma is elementary.
\begin{lemma}
 \label{lem:commuting_isotypes}
	If $I \subseteq J$, then for any $\sigma\in\Irr(K_{q,I})$ and $\tau\in\Irr(K_{q,J})$, the projections $p_\sigma$ and $p_\tau$ commute.
	In particular, $p_\sigma$ commutes with $p_\mu$ for any integral weight $\mu\in\bP$.
\end{lemma}

Note that the action $\rightharpoonup$ of $\sA(\hat{K_q})$ on $\sA(K_q)$ does not restrict to an action on the subspace $\sA(\sE_\mu)$.  But Lemma \ref{lem:commuting_isotypes} implies that the projection $p_\sigma\hit\,$ does restrict to a projection on $\sA(\sE_\mu)$.  In other words, it is still meaningful to talk of the decomposition of the section space $L^2(\sE_\mu)$ into isotypical components for the right regular action of the subgroup $K_{I,q}$.  This observation will allow us to define analogues of the algebras $\KK_I$ for the classical groups (Definition \ref{def:K-A-J}) which contain the longitudinal  pseudodifferential operators of negative order along the foliations $\sF_I$.

In the sequel, we write $p_\sigma$ for the operator $p_\sigma\hit\;$ on $L^2(\sE_\mu)$.
	
\begin{definition}
	Fix $q\in(0,1]$, let $\mu,\nu\in\bP$ and let $I\subseteq\bSigma$.  We define
  $\KK_I(L^2(\sE_\mu),L^2(\sE_\nu))$ to be the norm-closure of the set of bounded operators $T:L^2(\sE_\mu) \to L^2(\sE_\nu)$ satisfying
		\[
		 p_\sigma T = T p_\sigma = 0 \text{ for all but finitely many } \sigma \in \Irr(K_{I,q})
		\]
\end{definition}

This definition generalizes in an obvious way to operators between finite direct sums of line bundles $\sE = \bigoplus_i \sE_{\mu_i}$.   As in Chapter \ref{sec:flags}, we will often suppress the bundles in the notation, when they are clear from the context.  As a technical convenience, we can view the spaces $\KK_I(L^(\sE),L^2(\sE'))$ as the morphism spaces of a $C^*$-category $\KK_I$ with objects $L^2(\sE)$ for any $\sE$ direct sum of induced line bundles.

\begin{example}
  When $q=1$, the above definition of $\KK_I$ coincides with the algebra $\KK_I$ of Definition \ref{def:K-A-J}, that is, the hereditary $C^*$-algebra generated by the longitudinal smoothing operators $\Psi_I^{-\infty}$ along the fibres of the fibration $\sX \onto \sX_I$. 
\end{example}

The following theorem  is a quantum analogue of Theorem \ref{thm:foliations} for products of longitudinal pseudodifferential operators on classical flag manifolds.

\begin{theorem}[\cite{VoiYun:SUq3}]
 \label{thm:quantum_foliations}
 Let $K_q = \SU_q(n)$.  For any $I,J\subseteq\bSigma$, we have
 $\KK_I \KK_J \subseteq \KK_{I\cup J}$.
 Moreover, $\KK_\bSigma  = \KK$, the compact operators.
\end{theorem}

The proof of this theorem uses explicit calculations in terms of the Gelfand-Tsetlin bases of simple $\SU_q(n)$-modules, see \cite{VoiYun:SUq3}.
This type of calculation is not easy to generalize to other semisimple groups.  Although, one obviously expects that Theorem \ref{thm:quantum_foliations} should hold for all compact semisimple quantum groups, this remains an open problem.


\subsection{Principal series representations of $\SL_q(n,\CC)$}

As in Chapter \ref{sec:flags}, the above analysis is motivated by the equivariant index theory.   The Bernstein-Gelfand-Gelfand complex of Theorem \ref{thm:BGG} admits a $q$-deformation to a complex of unbounded intertwiners between non-unitary principal series representations of a quantized complex semisimple Lie group, see \cite{VoiYun:CQG}.   In this section, we will describe the normalized BGG complex, for which we will only need the \emph{base of principal series} representations, meaning those principal series with continuous parameter $\lambda=0$.    

We begin with some general notation.

\begin{definition}
Let $W \in M(\sA(K_q) \otimes \sAc(\hat{K}_q))$ denote the \emph{multiplicative unitary} of $K_q$, which is characterized by the following property.  For $f\in \sA(K_q)$ we write $[f] \in L^2(K_q)$ for the corresponding element in the GNS space.  We equip $L^2(K_q)$ the left actions $\lambda$ of $\sA(K_q)$ and $\hat\lambda$ of $\sAc(\hat{K}_q)$ defined by
\begin{align*}
 \lambda(a)[f] &= [af], && a,f\in\sA(K_q), \\
 \hat{\lambda}(x)[f] &= (\hat{S}(x),f_{(1)})[f_{(2)}], && x\in\sAc(\hat{K}_q), ~f\in\sA(K_q).
\end{align*}
Then $W$ satisfies 
\[
 (\lambda \otimes \hat\lambda)(W)\,([f] \otimes [g]) = [S^{-1}(g_{(1)})f] \otimes [g_{(2)}] ,
   \qquad\qquad \forall f,g\in\sA(K_q). 
\]
\end{definition}

\begin{definition}
The complex semisimple quantum group $G_q = \SL_q(n,\CC)$ is defined as the Drinfeld double of $K_q = \SU_q(n)$.  This means that it has function algebra
\[
 \sAc(G_q) = \sA(K_q) \bowtie \sAc(\hat{K}_q),
\]
which is equal to $ \sA(K_q) \otimes \sAc(\hat{K}_q)$ as an algebra, and equipped with the coalgebra operations
\begin{align*}
 \Delta_{G_q} &= (\id \otimes \mathrm{Flip} \otimes \id) (\id \otimes \ad(W) \otimes \id) (\Delta \otimes \hat\Delta), \\
 \epsilon_{G_q} &= \epsilon \otimes \hat\epsilon.
\end{align*}
\end{definition}

Note that $G_q$ is an algebraic quantum group in the sense of van Daele \cite{vanDaele:duality}.
We write $\sD(G_q)$ for the dual multiplier Hopf algebra, with operations dual to $\sAc(G_q)^\op$.

\begin{definition}
 Let $\mu\in\bP$.   The \emph{base of principal series representation} with parameter $\mu$ is the representation $\pi_{\mu,0}$ of $\sD(G_q)$ on the Hilbert space $L^2(\sE_\mu)$ which corresponds to the following $K_q$-Yetter-Drinfeld structure:
 \begin{itemize}
  \item coaction of $\sA(K_q)$ on $L^2(\sE_\mu)$ by the left regular corepresentation,
  \item action of $\sA(K_q)$ on $L^2(\sE_\mu)$ by the twisted adjoint action:
   \[
    a\cdot \xi = (K_{2\rho},a_{(2)}) \, a_{(1)} \xi S(a_{(3)})
   \]
 \end{itemize}
\end{definition}

The resulting representation of $\sD(G_q)$ is a $*$-representation, \emph{i.e.}\ corresponds to a unitary representation of $G_q$.  The following theorem summarizes some key properties of these representations, which are analogues of well-known results for the unitary principal series of the classical group $G$.

\begin{theorem}
\label{thm:G_q-intertwiners}~
 \begin{enumerate}
  \item The representations $\pi_{\mu,0}$ are all irreducible.
  \item Two representations $\pi_{\mu,0}$ and $\pi_{\mu',0}$ are unitarily equivalent if and only if $\mu'=w\mu$ for some $w\in W$.
  \item Suppose $w\to w'$ is an edge in the Bruhat graph associated to a simple root $\alpha$ and $\mu$ is a dominant weight.  Then the unitary intertwiner between $\pi_{-w\mu,0}$ and $\pi_{-w'\mu,0}$ is given by
  \[
    I_{w\to w'} :=
      \ph(E_\alpha^n)\hit\; : L^2(\sE_{-w\mu,0}) \to L^2(\sE_{-w'\mu,0}).
  \] 
  where $n$ is such that $w\mu-w'\mu = n\alpha$.
 \end{enumerate}
\end{theorem}

The first two results are essentially due to Joseph and Letzter \cite{JosLet:annihilators, Joseph:book}.   The third is quantum analogue of the formula \eqref{eq:intertwiners} when $\mu=\rho$.  For proofs of these, and many other related results, see \cite{VoiYun:CQG}.

\subsection{Algebras of longitudinal pseudodifferential operators, cont'd}

Morally speaking, the intertwiners of Theorem \ref{thm:G_q-intertwiners}(3) are longitudinal pseudodifferential operators of order $0$ along the fibres of $\sX_q \onto \sX_{\alpha,q}$.   Let us give a precise meaning to this. Following Definition \ref{def:K-A-J}, we define $\AA$ to be simultaneous multiplier category of all $\KK_I$ ($I\subseteq\bSigma$).  In other words, $A\in\AA_I(L^2(\sE,\sE'))$ if for any $I\subseteq\bSigma$ and any direct sum of line bundles $\sE''$ we have
\[
 A . \KK_I \subseteq \KK_I
 \qquad\text{and}\qquad
  \KK_I.A \subseteq \KK_I.
\]
We also put
\[
  \JJ_I = \KK_I \cap \AA.
\]
All of the operators we are interested in belong to $\AA$.

\begin{lemma}[\cite{VoiYun:SUq3}]
 \label{lem:quantum_A}
 Let $K_q=\SU_q(n)$.  For any $\mu\in\bP$, the following bounded operators belong to $\AA$:
 \begin{enumerate}
  \item[(i)] left or right multiplication by any $a\in C(\sE_\nu)$, as an operator $L^2(\sE_\mu) \to L^2(\sE_{\mu+\nu})$;
  \item[(ii)] the principal series representation $\pi_{\mu,0}(f)$ of any $f\in \sD(G_q)$ acting on $L^2(\sE_{\mu})$;
  \item[(iii)] the operator phases of the $\;\hit\;$ action of $E_\alpha$, $F_\alpha$ for any simple root $\alpha\in\bSigma$, viewed as operators $\ph(E_\alpha) : L^2(\sE_\mu) \to L^2(\sE_{\mu+\alpha})$ and $\ph(F_\alpha): L^2(\sE_\mu) \to L^2(\sE_{\mu-\alpha})$.
 \end{enumerate}
\end{lemma}

\begin{remark}
Parts \emph{(i)} and \emph{(ii)} of this lemma are straightforward, but the proof of part \emph{(iii)} relies on some rather delicate analysis using Gelfand-Tsetlin bases and $q$-orthogonal polynomials; see \cite{VoiYun:SUq3}.  Recall from Lemma \ref{lem:BBB_in_A} that this was one of the most difficult parts of the classical analysis, as well.
\end{remark}

Keeping in mind that the operators $\ph(E_\alpha^n) : L^2(\sE_\mu) \to L^2(\sE_{\mu+n\alpha})$ and $\ph(F_\alpha^n): L^2(\sE_\mu) \to L^2(\sE_{\mu-n\alpha})$ are to be thought of as order $0$ longitudinal pseudodifferential operators along the foliation $\sF_\alpha$, the following commutator relations should not be surprising.

\begin{lemma}
	\label{lem:phE-commutators}
 Let $\mu,\nu\in\bP$, $n \in \NN$.  For any section $a \in \sA(\sE_{\nu})$ we have 
\begin{align*}
 \ph(E^n_\alpha) a - a \ph(E^n_\alpha) & \in \JJ_\alpha(L^2(\sE_\mu), L^2(\sE_{\mu+\nu+n\alpha})), \\
 \ph(F^n_\alpha) a - a \ph(F^n_\alpha) &\in \JJ_\alpha(L^2(\sE_\mu), L^2(\sE_{\mu+\nu-n\alpha})).
\end{align*}
\end{lemma}

\subsection{The BGG class of the quantum flag manifold of $\SU_q(3)$}

Having obtained the key technical results---Theorem \ref{thm:quantum_foliations} and Lemma \ref{lem:quantum_A}---we can now translate the construction of the Kasparov module in Section \ref{sec:Kaparov_product} almost word-for-word to the quantum context. We will give an extremely brief summary.

Let $\alpha$, $\beta$ and $\rho=\alpha+\beta$ be the three positive roots of $\SU_q(3)$. 

\begin{lemma}
 The diagramme of \emph{normalized BGG operators}
 \begin{equation}
  \label{eq:q-BGG}
  \xymatrix@C=2ex@R=1ex{
   & L^2(\sE_{\alpha,0}) \ar[ddrr]^(0.8){\ph(F^2_\beta)} \ar[rr]^{T_1} \ar@{.}[dd]|{\oplus} && L^2(\sE_{2\alpha+\beta,0}) \ar[dr]^{\ph(E_\beta)} \ar@{.}@<2ex>[dd]|{\oplus}\\
   L^2(\sE_{0,0}) \ar[ur]^{\ph(E_\alpha)} \ar[dr]_{\ph(E_\beta)}
   &&&& L^2(\sE_{2\rho,0}), \\
   & L^2(\sE_{\beta,0}) \ar[uurr]_(0.8){\ph(F^2_\alpha)} \ar[rr]_{T_2} && L^2(\sE_{\alpha+2\beta,0}) \ar[ur]_{\ph(E_\alpha)}
 }
 \end{equation}
 with $T_1 = -\ph(E_\alpha^2)\ph(E_\beta)\ph(E_\alpha^*) $ and $T_2 = -\ph(E_\beta^2)\ph(E_\alpha)\ph(E_\beta^*)$,
 defines a complex modulo $\JJ_\alpha + \JJ_\beta$.  Moreover, the arrows which are pointed north-east (respectively, south-east, east) are invertible modulo $\JJ_\alpha$ (respectively, $\JJ_\beta$, $\JJ_\alpha+\JJ_\beta$).
\end{lemma}

\begin{proof}
 By Hopf-Galois theory, we can find sections $f_1,\ldots,f_n \in \sA(\sE_{-\rho})$ such that $\sum_{i=1}^n f_i^*f_i = 1 \in \sA(\sX)$.  From Lemma \ref{lem:phE-commutators} we obtain that
 \[
  \left(\sum_i f_i^* \ph(E_\alpha) f_i\right) - \ph(E_\alpha)  \quad \in\JJ_\alpha(L^2(\sE_\mu),L^2(\sE_{\mu+\alpha}),
 \]
 for any $\mu\in\bP$.  The result then follows, as in the classical case of Lemma \ref{lem:shift_trick}, by comparing the diagram of normalized BGG operators \eqref{eq:q-BGG} with the diagram of intertwiners $I_{w\to w'}$ from Theorem \ref{thm:G_q-intertwiners}.
\end{proof}

We now have all the elements that we used in the construction of the BGG class for classical $\SL(3,\CC)$.  We can repeat the application of the Kasparov Technical Theorem to obtain operators $N_\alpha,N_\beta$ analogous to those in Equation \eqref{eq:po1}, and then modify the normalized BGG complex as in Figure \ref{fig:KTT}.  
The result is an Kasparov $K$-homology class $\BGG$ for $\sX_q$ which is equivariant with respect to the Drinfeld double $\SL_q(3,\CC)$ of $\SU_q(3)$.  The details can be found in \cite{VoiYun:SUq3}.

\begin{theorem} \textup{(\cite{VoiYun:SUq3})}
	The above construction yields an $\SL_q(3,\CC)$-equivariant $K$-homology class for $\sX_q$, which we denote by $\BGG \in KK^{\SL_q(3,\CC)}(C(\sX_q), \CC)$.  Its image under the forgetful map to $KK^{\SU_q(3)}(\CC ,\CC)$ is the class of the trivial representation.
\end{theorem}

Let us finish with two immediate applications of the BGG class.

\medskip

The first application is the equivariant Poincar\'e duality of the flag manifold $\sX_q$.  The notion of $KK$-theoretic Poincar\'e duality was introduced by Kasparov, but its generalization to quantum group equivariant $KK$-theory requires some considerable technical refinement.  This is due to the observation of Nest and Voigt \cite{NesVoi} that the Kasparov product in $KK^{K_q}$ requires the use of the braided tensor product $\boxtimes$.  Therefore, even if we want to prove only $K_q$-equivariant Poincar\'e duality, we are obliged to define a Poincar\'e duality class in $KK^{G_q}$, where $G_q$ is the Drinfeld double of $K_q$.

Nest and Voigt \cite{NesVoi} showed how to prove $\SU_q(2)$-equivariant Poincar\'e duality for the Podle\'{s} sphere, using an $\SL_q(2,\CC)$-equivariant duality class.  The same argument goes through for the full flag variety of $\SU_q(3)$ using the class $\BGG$, and yields the following.

\begin{theorem}
  The quantum flag manifold $ \sX_q $ is $ \SU_q(3) $-equivariantly Poincar\'e dual to itself. That is, there is a natural isomorphism
  $$
  KK^{\SL_q(3,\CC)}_*(C(\sX_q) \boxtimes A, B) \cong KK^{\SL_q(3,\CC)}_*(A, C(\sX_q) \boxtimes B)
  $$
  for all $ \SL_q(3,\CC) $-$ C^* $-algebras $ A $ and $ B $.
\end{theorem}

\medskip

The second application is to the Baum-Connes Conjecture.  Once again, making sense of the Baum-Connes Conjecture for discrete quantum groups requires some serious technical machinery, this time developed by Meyer and Nest \cite{MeyNes:localisation, MeyNes:coactions}.  The Baum-Connes Conjecture for the discrete dual of $\SU_q(2)$ was proven by Voigt in \cite{Voigt:free_orthogonal}.  An analogous argument using the class $\BGG$ yields the following.

\begin{theorem}
	The discrete quantum dual of $\SU_q(3)$ with $q \in (0,1]$ satisfies the Baum-Connes Conjecture with trivial coefficients.
\end{theorem}

%

\bibliographystyle{alpha}
\bibliography{refs}

\def\polhk#1{\setbox0=\hbox{#1}{\ooalign{\hidewidth
  \lower1.5ex\hbox{`}\hidewidth\crcr\unhbox0}}}
\begin{thebibliography}{NWX99}

\bibitem[BE89]{BasEas}
Robert~J. Baston and Michael~G. Eastwood.
\newblock {\em The {P}enrose transform}.
\newblock Oxford Mathematical Monographs. The Clarendon Press, Oxford
  University Press, New York, 1989.
\newblock Its interaction with representation theory, Oxford Science
  Publications.

\bibitem[BG88]{BeaGre}
Richard Beals and Peter Greiner.
\newblock {\em Calculus on {H}eisenberg manifolds}, volume 119 of {\em Annals
  of Mathematics Studies}.
\newblock Princeton University Press, Princeton, NJ, 1988.

\bibitem[BGG75]{BerGelGel}
I.~N. Bernstein, I.~M. Gelfand, and S.~I. Gelfand.
\newblock Differential operators on the base affine space and a study of
  {${\mathfrak{g}}$}-modules.
\newblock pages 21--64, 1975.

\bibitem[BN03]{BenNis}
Moulay-Tahar Benameur and Victor Nistor.
\newblock Homology of algebras of families of pseudodifferential operators.
\newblock {\em J. Funct. Anal.}, 205(1):1--36, 2003.

\bibitem[Boh15]{Bohlen:bisingular}
Karsten Bohlen.
\newblock The {$K$}-theory of bisingular pseudodifferential algebras.
\newblock {\em J. Pseudo-Differ. Oper. Appl.}, 6(3):361--382, 2015.

\bibitem[CF03]{CraFer}
Marius Crainic and Rui~Loja Fernandes.
\newblock Integrability of {L}ie brackets.
\newblock {\em Ann. of Math. (2)}, 157(2):575--620, 2003.

\bibitem[CM95]{ConMos:local_index_formula}
A.~Connes and H.~Moscovici.
\newblock The local index formula in noncommutative geometry.
\newblock {\em Geom. Funct. Anal.}, 5(2):174--243, 1995.

\bibitem[CM08]{ConMos:twisted}
Alain Connes and Henri Moscovici.
\newblock Type {III} and spectral triples.
\newblock In {\em Traces in number theory, geometry and quantum fields},
  Aspects Math., E38, pages 57--71. Friedr. Vieweg, Wiesbaden, 2008.

\bibitem[Con94]{Connes:NCG}
Alain Connes.
\newblock {\em Noncommutative geometry}.
\newblock Academic Press, Inc., San Diego, CA, 1994.

\bibitem[CP]{ChoPon}
Woocheol Choi and Raphael Ponge.
\newblock Privileged coordinates and tangent groupoid for carnot manifolds.
\newblock Preprint. \url{http://arxiv.org/abs/1510.05851}.

\bibitem[{\v{C}}S09]{CapSlo}
Andreas {\v{C}}ap and Jan Slov{\'a}k.
\newblock {\em Parabolic geometries. {I}}, volume 154 of {\em Mathematical
  Surveys and Monographs}.
\newblock American Mathematical Society, Providence, RI, 2009.
\newblock Background and general theory.

\bibitem[{\v{C}}SS01]{CapSloSou:BGG}
Andreas {\v{C}}ap, Jan Slov{\'{a}}k, and Vladim{\'{\i}}r Sou{\v{c}}ek.
\newblock Bernstein-{G}elfand-{G}elfand sequences.
\newblock {\em Ann. of Math. (2)}, 154(1):97--113, 2001.

\bibitem[DD10]{DAnDab}
Francesco D'Andrea and Ludwik D{\polhk{a}}browski.
\newblock Dirac operators on quantum projective spaces.
\newblock {\em Comm. Math. Phys.}, 295(3):731--790, 2010.

\bibitem[Deb01]{Debord:Integration}
Claire Debord.
\newblock Local integration of {L}ie algebroids.
\newblock In {\em Lie algebroids and related topics in differential geometry
  ({W}arsaw, 2000)}, volume~54 of {\em Banach Center Publ.}, pages 21--33.
  Polish Acad. Sci., Warsaw, 2001.

\bibitem[DH17]{DavHal:BGG}
Shantanu Dave and Stefan Haller.
\newblock Graded hypoellipticity of bgg sequences.
\newblock Preprint. \url{http://arxiv.org/abs/1705.01659}, 2017.

\bibitem[DS14]{DebSka:Rx-action}
Claire Debord and Georges Skandalis.
\newblock Adiabatic groupoid, crossed product by {$\Bbb{R}_+^\ast$} and
  pseudodifferential calculus.
\newblock {\em Adv. Math.}, 257:66--91, 2014.

\bibitem[Duf75]{Duflo:representations}
Michel Duflo.
\newblock Repr\'esentations irr\'eductibles des groupes semi-simples complexes.
\newblock pages 26--88. Lecture Notes in Math., Vol. 497, 1975.

\bibitem[EM]{EpsMel:book}
Charles Epstein and Richard Melrose.
\newblock The heisenberg algebra, index theory and homology.
\newblock Unpublished book, \url{http://www-math.mit.edu/~rbm/book.html}.

\bibitem[EMM91]{EpsMelMen}
C.~L. Epstein, R.~B. Melrose, and G.~A. Mendoza.
\newblock Resolvent of the {L}aplacian on strictly pseudoconvex domains.
\newblock {\em Acta Math.}, 167(1-2):1--106, 1991.

\bibitem[FS74]{FolSte:estimates}
G.~B. Folland and E.~M. Stein.
\newblock Estimates for the {$\bar \partial _{b}$} complex and analysis on the
  {H}eisenberg group.
\newblock {\em Comm. Pure Appl. Math.}, 27:429--522, 1974.

\bibitem[FS82]{FolSte:Hardy}
G.~B. Folland and Elias~M. Stein.
\newblock {\em Hardy spaces on homogeneous groups}, volume~28 of {\em
  Mathematical Notes}.
\newblock Princeton University Press, Princeton, N.J.; University of Tokyo
  Press, Tokyo, 1982.

\bibitem[Hig04]{Higson:local_index_formula}
Nigel Higson.
\newblock The local index formula in noncommutative geometry.
\newblock In {\em Contemporary developments in algebraic {$K$}-theory}, ICTP
  Lect. Notes, XV, pages 443--536. Abdus Salam Int. Cent. Theoret. Phys.,
  Trieste, 2004.

\bibitem[HN79]{HelNou}
B.~Helffer and J.~Nourrigat.
\newblock Caracterisation des op\'erateurs hypoelliptiques homog\`enes
  invariants \`a gauche sur un groupe de {L}ie nilpotent gradu\'e.
\newblock {\em Comm. Partial Differential Equations}, 4(8):899--958, 1979.

\bibitem[Hum08]{Humphreys:O}
James~E. Humphreys.
\newblock {\em Representations of semisimple {L}ie algebras in the {BGG}
  category {$\mathcal{O}$}}, volume~94 of {\em Graduate Studies in
  Mathematics}.
\newblock American Mathematical Society, Providence, RI, 2008.

\bibitem[JK95]{JulKas}
Pierre Julg and Gennadi Kasparov.
\newblock Operator {$K$}-theory for the group {${\rm SU}(n,1)$}.
\newblock {\em J. Reine Angew. Math.}, 463:99--152, 1995.

\bibitem[JL95]{JosLet:annihilators}
Anthony Joseph and Gail Letzter.
\newblock Verma module annihilators for quantized enveloping algebras.
\newblock {\em Ann. Sci. \'Ecole Norm. Sup. (4)}, 28(4):493--526, 1995.

\bibitem[Jos95]{Joseph:book}
Anthony Joseph.
\newblock {\em Quantum groups and their primitive ideals}, volume~29 of {\em
  Ergebnisse der Mathematik und ihrer Grenzgebiete (3) [Results in Mathematics
  and Related Areas (3)]}.
\newblock Springer-Verlag, Berlin, 1995.

\bibitem[Jul02]{Julg:Spn1}
Pierre Julg.
\newblock La conjecture de {B}aum-{C}onnes \`a coefficients pour le groupe
  {${\rm Sp}(n,1)$}.
\newblock {\em C. R. Math. Acad. Sci. Paris}, 334(7):533--538, 2002.

\bibitem[Kas84]{Kasparov:Lorentz}
G.~G. Kasparov.
\newblock Lorentz groups: {$K$}-theory of unitary representations and crossed
  products.
\newblock {\em Dokl. Akad. Nauk SSSR}, 275(3):541--545, 1984.
\newblock English translation: Soviet Math. Dokl. 29 (1984), no. 2, 256--260.

\bibitem[Kna01]{Knapp:representation_theory}
Anthony~W. Knapp.
\newblock {\em Representation theory of semisimple groups}.
\newblock Princeton Landmarks in Mathematics. Princeton University Press,
  Princeton, NJ, 2001.
\newblock An overview based on examples, Reprint of the 1986 original.

\bibitem[Kr{\"a}04]{Krahmer:Dirac}
Ulrich Kr{\"a}hmer.
\newblock Dirac operators on quantum flag manifolds.
\newblock {\em Lett. Math. Phys.}, 67(1):49--59, 2004.

\bibitem[KTS15]{KraTuc}
Ulrich Kr{\"a}hmer and Matthew Tucker-Simmons.
\newblock On the {D}olbeault-{D}irac operator of quantized symmetric spaces.
\newblock {\em Trans. London Math. Soc.}, 2(1):33--56, 2015.

\bibitem[LMV17]{LesManVas}
Jean-Marie Lescure, Dominique Manchon, and St\'ephane Vassout.
\newblock About the convolution of distributions on groupoids.
\newblock {\em J. Noncommut. Geom.}, 11(2):757--789, 2017.

\bibitem[Mel82]{Melin:preprint}
Anders Melin.
\newblock Lie filtrations and pseudo-differential operators.
\newblock Preprint, 1982.

\bibitem[MN06]{MeyNes:localisation}
Ralf Meyer and Ryszard Nest.
\newblock The {B}aum-{C}onnes conjecture via localisation of categories.
\newblock {\em Topology}, 45(2):209--259, 2006.

\bibitem[MN07]{MeyNes:coactions}
Ralf Meyer and Ryszard Nest.
\newblock An analogue of the {B}aum-{C}onnes isomorphism for coactions of
  compact groups.
\newblock {\em Math. Scand.}, 100(2):301--316, 2007.

\bibitem[Mos10]{Moscovici:twisted}
Henri Moscovici.
\newblock Local index formula and twisted spectral triples.
\newblock In {\em Quanta of maths}, volume~11 of {\em Clay Math. Proc.}, pages
  465--500. Amer. Math. Soc., Providence, RI, 2010.

\bibitem[MY]{MatYun:PsiDOs}
Marco Matassa and Robert Yuncken.
\newblock Regularity of twisted spectral triples and pseudodifferential
  calculi.
\newblock {\em J. Noncommut. Geom.}
\newblock To appear. Preprint at \url{https://arxiv.org/abs/1705.04178}.

\bibitem[Nis00]{Nistor:integration}
Victor Nistor.
\newblock Groupoids and the integration of {L}ie algebroids.
\newblock {\em J. Math. Soc. Japan}, 52(4):847--868, 2000.

\bibitem[NT05]{NesTus:local_index_formula}
Sergey Neshveyev and Lars Tuset.
\newblock A local index formula for the quantum sphere.
\newblock {\em Comm. Math. Phys.}, 254(2):323--341, 2005.

\bibitem[NT10]{NesTus:isospectral}
Sergey Neshveyev and Lars Tuset.
\newblock The {D}irac operator on compact quantum groups.
\newblock {\em J. Reine Angew. Math.}, 641:1--20, 2010.

\bibitem[NV10]{NesVoi}
Ryszard Nest and Christian Voigt.
\newblock Equivariant {P}oincar\'e duality for quantum group actions.
\newblock {\em J. Funct. Anal.}, 258(5):1466--1503, 2010.

\bibitem[NWX99]{NisWeiXu}
Victor Nistor, Alan Weinstein, and Ping Xu.
\newblock Pseudodifferential operators on differential groupoids.
\newblock {\em Pacific J. Math.}, 189(1):117--152, 1999.

\bibitem[Pon06]{Ponge:groupoid}
Rapha{\"e}l Ponge.
\newblock The tangent groupoid of a {H}eisenberg manifold.
\newblock {\em Pacific J. Math.}, 227(1):151--175, 2006.

\bibitem[Pon08]{Ponge:Heisenberg}
Rapha{\"e}l~S. Ponge.
\newblock Heisenberg calculus and spectral theory of hypoelliptic operators on
  {H}eisenberg manifolds.
\newblock {\em Mem. Amer. Math. Soc.}, 194(906):viii+ 134, 2008.

\bibitem[Pus11]{Puschnigg:T}
Michael Puschnigg.
\newblock Finitely summable {F}redholm modules over higher rank groups and
  lattices.
\newblock {\em J. K-Theory}, 8(2):223--239, 2011.

\bibitem[Roc78]{Rockland}
Charles Rockland.
\newblock Hypoellipticity on the {H}eisenberg group-representation-theoretic
  criteria.
\newblock {\em Trans. Amer. Math. Soc.}, 240:1--52, 1978.

\bibitem[Rod75]{Rodino:bisingular}
Luigi Rodino.
\newblock A class of pseudo differential operators on the product of two
  manifolds and applications.
\newblock {\em Ann. Scuola Norm. Sup. Pisa Cl. Sci. (4)}, 2(2):287--302, 1975.

\bibitem[Tay]{Taylor:microlocal}
Michael Taylor.
\newblock Noncommutative microlocal analysis, part {I} (revised version).
\newblock \url{http://www.unc.edu/math/Faculty/met/NCMLMS.pdf}.

\bibitem[Tr{\`e}67]{Treves:TVS}
Fran{\c{c}}ois Tr{\`e}ves.
\newblock {\em Topological vector spaces, distributions and kernels}.
\newblock Academic Press, New York-London, 1967.

\bibitem[Uuy11]{Uuye:PsiDOs}
Otgonbayar Uuye.
\newblock Pseudo-differential operators and regularity of spectral triples.
\newblock In {\em Perspectives on noncommutative geometry}, volume~61 of {\em
  Fields Inst. Commun.}, pages 153--163. Amer. Math. Soc., Providence, RI,
  2011.

\bibitem[Vae05]{Vaes:induction}
Stefaan Vaes.
\newblock A new approach to induction and imprimitivity results.
\newblock {\em J. Funct. Anal.}, 229(2):317--374, 2005.

\bibitem[VD98]{vanDaele:duality}
A.~Van~Daele.
\newblock An algebraic framework for group duality.
\newblock {\em Adv. Math.}, 140(2):323--366, 1998.

\bibitem[vE05]{VanErp:thesis}
Erik van Erp.
\newblock {\em The {A}tiyah-{S}inger index formula for subelliptic operators on
  contact manifolds}.
\newblock ProQuest LLC, Ann Arbor, MI, 2005.
\newblock Thesis (Ph.D.)--The Pennsylvania State University.

\bibitem[vEY]{VanYun:PsiDOs}
Erik van Erp and Robert Yuncken.
\newblock A groupoid approach to pseudodifferential operators.
\newblock {\em J. Reine Agnew. Math}.
\newblock To appear. \url{https://doi.org/10.1515/crelle-2017-0035}.

\bibitem[vEY17]{VanYun:groupoid}
Erik van Erp and Robert Yuncken.
\newblock On the tangent groupoid of a filtered manifold.
\newblock {\em Bull. Lond. Math. Soc.}, 49(6):1000--1012, 2017.

\bibitem[Voi11]{Voigt:free_orthogonal}
Christian Voigt.
\newblock The {B}aum-{C}onnes conjecture for free orthogonal quantum groups.
\newblock {\em Adv. Math.}, 227(5):1873--1913, 2011.

\bibitem[VY15]{VoiYun:SUq3}
Christian Voigt and Robert Yuncken.
\newblock Equivariant {F}redholm modules for the full quantum flag manifold of
  {${\rm SU}_q(3)$}.
\newblock {\em Doc. Math.}, 20:433--490, 2015.

\bibitem[VY17]{VoiYun:CQG}
Christian Voigt and Robert Yuncken.
\newblock Complex semisimple quantum groups and representation theory.
\newblock Preprint. \url{https://arxiv.org/abs/1705.05661} (226 pages), 2017.

\bibitem[Yun10]{Yuncken:products_of_PsiDOs}
Robert Yuncken.
\newblock Products of longitudinal pseudodifferential operators on flag
  varieties.
\newblock {\em J. Funct. Anal.}, 258(4):1140--1166, 2010.

\bibitem[Yun11]{Yuncken:BGG}
Robert Yuncken.
\newblock The {B}ernstein-{G}elfand-{G}elfand complex and {K}asparov theory for
  {${\rm SL}(3,\Bbb C)$}.
\newblock {\em Adv. Math.}, 226(2):1474--1512, 2011.

\bibitem[Yun13]{Yuncken:foliations}
Robert Yuncken.
\newblock Foliation {$C^*$}-algebras on multiply fibred manifolds.
\newblock {\em J. Funct. Anal.}, 265(9):1829--1839, 2013.

\end{thebibliography}

\newpage 
\thispagestyle{empty}~
\newpage 
\thispagestyle{empty}~

\end{document}